\documentclass[11pt]{amsart} 

\usepackage{amsmath,amssymb}
\usepackage[mathscr]{eucal}
\usepackage[all]{xy}

\newcommand{\epi}{\twoheadrightarrow}
\newcommand{\iso}{\buildrel{\sim}\over{\longrightarrow}}
\newcommand{\mono}{\hookrightarrow}

\newcommand{\limto}{{\displaystyle\lim_{\longrightarrow}}}
\newcommand{\rightlim}{\mathop{\limto}}

\newcommand{\leftlim}{\mathop{\displaystyle\lim_{\longleftarrow}}}
\newcommand{\limfromn}{\leftlim\limits_{\raise3pt\hbox{$n$}}}
\newcommand{\limton}{\rightlim\limits_{\raise3pt\hbox{$n$}}}

\newcommand{\rightlimit}[1]{\mathop{\lim\limits_{\longrightarrow}}\limits%
                   _{\raise3pt\hbox{$\scriptstyle #1$}}}
  
\newcommand{\leftlimit}[1]{\mathop{\lim\limits_{\longleftarrow}}\limits%
                   _{\raise3pt\hbox{$\scriptstyle #1$}}}

\newcommand{\BA}{{\mathbb A}}
\newcommand{\BC}{{\mathbb C}}

\newcommand{\BN}{{\mathbb N}}
\newcommand{\BP}{{\mathbb P}}
\newcommand{\BQ}{{\mathbb Q}}
\newcommand{\BR}{{\mathbb R}}
\newcommand{\BZ}{{\mathbb Z}}

\newcommand{\cA}{{\mathcal A}}
\newcommand{\cB}{{\mathcal B}}
\newcommand{\cC}{{\mathcal C}}
\newcommand{\cD}{{\mathcal D}}
\newcommand{\cG}{{\mathcal G}}

\newcommand{\cL}{{\mathcal L}}
\newcommand{\cK}{{\mathcal K}}
\newcommand{\cM}{{\mathcal M}}
\newcommand{\cO}{{\mathcal O}}
\newcommand{\cP}{{\mathcal P}}
\newcommand{\cR}{{\mathcal R}}
\newcommand{\cS}{{\mathcal S}}
\newcommand{\cT}{{\mathcal T}}

\newcommand{\IS}{Tate-smooth}
\newcommand{\MLT}{Mittag-Leffler-Tate }

\newcommand{\bd}{{\big\downarrow}}

\newcommand{\hatimes}{{\widehat{\otimes}}}

\newcommand{\coflat}{coprojective }
\newcommand{\allmost}{almost }
\newcommand{\ALMOST}{2-almost }

\DeclareMathOperator{\WP}{{ap}}

\DeclareMathOperator{\Aut}{{Aut}}
\DeclareMathOperator{\AAut}{{\mathcal A\it ut}}
\DeclareMathOperator{\Calk}{{Calk}}
\DeclareMathOperator{\Card}{{Card}}
\DeclareMathOperator{\can}{{can}}
\DeclareMathOperator{\Cl}{{Cl}}
\DeclareMathOperator{\comp}{{comp}}

\DeclareMathOperator{\Cone}{{Cone}}
\DeclareMathOperator{\Char}{{char}}
\DeclareMathOperator{\Dim}{{Dim}}
\DeclareMathOperator{\discr}{{discr}}
\DeclareMathOperator{\Det}{{Det}}
\DeclareMathOperator{\End}{{End}}
\DeclareMathOperator{\et}{{et}}
\DeclareMathOperator{\Ext}{{Ext}}
\DeclareMathOperator{\Funct}{{Funct}}
\DeclareMathOperator{\Gal}{{Gal}}
\DeclareMathOperator{\GR}{{\,\cG\cR}}
\DeclareMathOperator{\Gras}{{Gras}}
\DeclareMathOperator{\Gun}{{Gun}}
\DeclareMathOperator{\hatAut}{{\widehat{\mathcal A}ut}}

\DeclareMathOperator{\Hom}{{Hom}}
\DeclareMathOperator{\id}{{id}}
\DeclareMathOperator{\im}{{Im}}

\DeclareMathOperator{\K}{{K}}

\DeclareMathOperator{\Kar}{{Kar}}
\DeclareMathOperator{\Ker}{{Ker}}
\DeclareMathOperator{\Lie}{{Lie}}

\DeclareMathOperator{\Nis}{{Nis}}

\DeclareMathOperator{\Ob}{{Ob}}
\DeclareMathOperator{\Pic}{{Pic}}

\DeclareMathOperator{\Pregun}{{Pre-gun}}
\DeclareMathOperator{\RR}{{R}}
\DeclareMathOperator{\rank}{{rank}}
\DeclareMathOperator{\red}{{red}}
\DeclareMathOperator{\res}{{res}}
\DeclareMathOperator{\Spec}{{Spec}}

\DeclareMathOperator{\Top}{{top}}

\DeclareMathOperator{\Ttop}{{-top}}
\DeclareMathOperator{\TTop}{{topol}}

\DeclareMathOperator{\Vect}{{Vect}}

\newcommand{\RGamma}{(R\Gamma)_{\TTop}}
\newcommand{\bcdot}{{\textstyle\cdot}}

\DeclareMathOperator{\all}{{all}}

\newcommand{\BG}{{\mathbb G}}

\newcommand{\fU}{{\mathfrak{U}}}

\newcommand{\cF}{{\mathcal F}}
\newcommand{\Call}{\cC^{\all}}

\newtheorem{pro}{Proposition}[section]
\newtheorem{lm}[pro]{Lemma}
\newtheorem{theor}[pro]{Theorem}

\numberwithin{equation}{section}

\begin{document}

\thanks{Partially supported by NSF grants DMS-0100108 and
DMS-0401164.}
\address{Dept. of Math., Univ. of Chicago, 5734 University 
Ave., Chicago, IL 60637, USA}

\title{Infinite-dimensional vector bundles in algebraic
geometry (an introduction)}
\author{Vladimir Drinfeld}

\email{drinfeld@math.uchicago.edu}

\dedicatory{To Izrail Moiseevich Gelfand with deepest
gratitude and admiration}

\maketitle

\tableofcontents

{{{\leftskip=2.00in {
{{\noindent \it  According to the Lamb conjecture, the
key to the future development of Quantum Field  Theory is
probably buried in some forgotten paper published in the 30's.
Attempts to follow up this conjecture, however, will probably
be  unsuccessful because of the Peierls-Jensen paradox;
namely, that even if one finds the right paper, the point
will probably be missed until it is found independently and
accidentally by experiment.}

\smallskip 

\noindent  The Future of Field Theory, by Pure Imaginary
Observer \cite{PI}. 
}}\par}

}}

\bigskip

\section{Introduction}
\subsection{Subject of the article}
The goal of this work is to show that there is a reasonable
algebro-geometric notion of vector bundle
with infinite-dimensional locally linearly compact fibers and
that these objects appear ``in nature". Our approach is based
on some results and ideas discovered in algebra in 1958--1972
by H.~Bass, L.~Gruson, I.~Kaplansky, M.~Karoubi, and
M.~Raynaud.

This article contains definitions and formulations of the main
theorems, but practically no proofs. A detailed exposition
will appear in \cite{Dr}.

\subsection{Conventions} We use the words {\it ``$S$-family
of vector spaces''\,} as shorthand for ``vector bundle on a
scheme $S$'' and {\it ``Tate space"\,} as shorthand for
``locally linearly compact vector space".

\subsection{Overview of the results and structure of the
article}
\subsubsection{General theory}
In \S\ref{discr} we recall the Raynaud-Gruson theorem on the
local nature of projectivity, which shows that there is a 
good notion of family of discrete infinite-dimensional vector
spaces.

In \S\ref{Tatesection} we introduce the notion of {\it ``Tate
module"\,} over an arbitrary ring $R$ and show that if $R$ is
commutative one thus gets a reasonable notion of $S$-family
of Tate spaces, $S=\Spec R$. One has to take in account that
$K_0$ of the additive category of Tate $R$-modules may be
nontrivial. In fact, it equals $K_{-1}(R)$. We show that
$K_{-1}(R)=0$ if $R$ is Henselian. We give a proof of this
fact because it explains the fundamental role of the
Nisnevich topology in this work. We discuss the notions of
dimension torsor and determinant gerbe of a family of Tate
spaces.

At least technically, the theory of Tate $R$-modules is based
on the notion of {\it\allmost projective\,} module, which is
introduced in \S\ref{Almostsect}. Roughly speaking, a module
is \allmost projective if it is projective up to finitely
generated modules. Unlike Tate modules, \allmost projective
modules are discrete. Any Tate module can be represented as
the projective limit of a filtering projective system of
\allmost projective modules with surjective transition maps. 

\S\ref{dim&det} is devoted to the canonical central extension
of the automorphism groups of \allmost projective and
Tate $R$-modules. In \ref{Beilvision} we discuss an
interesting (though slightly vague) picture, which I learned
from A.~Beilinson.

\subsubsection{Application to the space of formal loops}
In \S\ref{refinedsection} we define a class of 
{\it Tate-smooth\,} ind-schemes (morally, these are smooth
infinite-dimensional algebraic manifolds modeled on Tate
spaces). According to Theorem \ref{ind-smooth}, the 
ind-scheme of formal loops of a smooth affine manifold $Y$
over the local field $k((t))$ is Tate-smooth over $k$. This
is one of our main results. In \ref{refined} we use it to
define a ``refined"  version of the motivic integral of a
differential form on $Y$ with no zeros over the intersection
of $Y$ with a polydisk. Unlike the usual motivic integral,
the ``refined" one is an object of a triangulated category
rather than an element of a group.

\subsubsection{Application to vector bundles on a manifold
with punctures} In \S\ref{polydisksection} we first show that
\allmost projective and Tate modules appear naturally in the
study of the cohomology of a family of finite-dimensional
vector bundles on a punctured smooth manifold. Then we
briefly explain how the canonical central extension that comes
from this cohomology allows (in the case of $GL(n)$-bundles)
to interpret the ``Uhlenbeck compactification" constructed in
\cite{FGK,BFG} as the fine moduli space of a certain type of
generalized vector bundles on $\BP^2$ (we call them {\it
gundles\,}). In fact, the application to the ``Uhlenbeck
compactification" was one of the main motivations of this
work.

\subsection{Acknowledgements}
I thank A.~Beilinson, D.~Gaitsgory,
V.~Ginzburg, D.~Hirschfeldt, M.~Kapranov, D.~Kazhdan, 
A.~Suslin, and C.~Weibel for stimulating discussions.

Speaking at the ``Unity of Mathematics'' conference (Harvard,
2003) was a great honor and pleasure for me, and an
important stimulus to write this article. I am very grateful
to the organizers of the conference.

\section{Families of discrete infinite-dimensional vector
spaces (after Kaplansky, Raynaud and Gruson).}
\label{discr}

Is there a reasonable notion of not necessarily
finite-dimensional vector bundle on a scheme? We know due to
Serre \cite{S} that a finite-dimensional vector bundle on
an affine scheme $\Spec R$ is the same as a finitely
generated projective $R$-module. So it is natural to give the
following definition.

\medskip

\noindent {\bf Definition.}
A vector bundle on a scheme $X$ is a quasicoherent sheaf of
$\mathcal{O}_X$-modules $\mathcal{F}$ such that for every
open affine subset
$\Spec R\subset X$ the $R$-module
$H^0(\operatorname{Spec}R,\mathcal{F})$ is projective.

\medskip

{\it Key question:} is this a local notion?
More precisely, the question is as follows. Let
$\operatorname{Spec}R=\bigcup\limits_iU_i$, $U_i=\operatorname{Spec}R_i$.
Let $M$ be a (not necessarily finitely generated) $R$-module
such that $M\otimes_RR_i$ is  projective for all $i$. Does it
follow that $M$ is projective?

The question is difficult: the arguments used in the case that
$M$ is finitely generated fail for modules of infinite type.
Nevertheless Grothendieck (\cite{Gr}, Remark 9.5.8)
conjectured that the answer is positive. This was proved by
Raynaud and Gruson \cite{RG} (in Ch.~1 for countably generated
modules and in Ch.~2 for arbitrary ones). Moreover, they
proved the following theorem, which says that projectivity is
a local property for the fpqc topology (not only for Zariski).

\begin{theor} \label{Mit.1}
Let $M$ be a module over a commutative ring $R$ and $R'$ be
a flat commutative $R$-algebra such that the morphism 
$\Spec R'\to\Spec R$ is surjective. If $R'\otimes_RM$ is
projective then $M$ is. 
\end{theor}

In fact, they derived it as an easy
corollary of the following remarkable and nontrivial theorem
due to Kaplansky \cite{Ka} and Raynaud-Gruson \cite{RG},
which explains what projectivity really is. Theorem
\ref{Mit.1} follows from the fact that for commutative rings
properties (a)-(c) below are local.

\begin{theor} \label{Mit.2}
Let $R$ be a (not necessarily commutative) ring.
An $R$-module $M$ is projective if and only if the following
properties hold:

(a) $M$ is flat;

(b) $M$ is a direct sum of countably generated modules;

(c) $M$ is a Mittag-Leffler module.
\end{theor}

The fact that a projective module can be represented as a
direct sum of countably generated ones was proved by
Kaplansky \cite{Ka}. 

\medskip

The remaining part of Theorem \ref{Mit.2} is due to Raynaud
and Gruson \cite{RG}. The key notion of Mittag-Leffler module
was introduced in \cite{RG}. Here I prefer only to explain
what a flat Mittag-Leffler module is. By the Govorov--Lazard
lemma \cite{Gov,Laz}, a flat $R$-module $M$ can be represented
as the inductive limit of a directed family of finitely
generated projective modules $P_i$. According to \cite{RG}, in
this situation $M$ is Mittag-Leffler if and only if the
projective system formed by the dual (right) $R$-modules
$P_i^*:=\Hom_R(P_i,R)$ satisfies the Mittag-Leffler
condition: for every $i$ there exists  $j\ge i$ such that 
$\im (P^*_j\to P_i^*)=\im (P^*_k\to P_i^*)$ for all $k\ge j$.

\medskip

\noindent {\bf Remarks.} 
(i) One gets a sightly different definition of not
necessarily finite-dimensional vector bundle on a scheme
if one replaces projectivity by the property of being a
flat Mittag-Leffler module. The product of infinitely many
copies of $\BZ$ is an example of a flat Mittag-Leffler 
$\BZ$-module which is not projective (it is due to Baer,
see p.48 and p.82 of \cite{Ka2}).
 Unlike projectivity,
the property of $M$ being a flat Mittag-Leffler module is
a first-order property (in the sense of mathematical logic)
of $R^{(\BN )}\otimes_RM$ viewed as a module over 
$\End_RR^{(\BN )}$ (here $R^{(\BN )}$ is the right
$R$-module freely generated by $\BN$). Let me also mention
that one does not need AC (the axiom of choice) to prove
that a vector  space over a field is a flat Mittag-Leffler
module, but in set theory without AC one cannot prove that
$\BR$ is a direct summand of a free
$\BQ$-module\footnote{Without AC it is not true that any
free module $F$ is projective, i.e., every epimorphism $M\to
F$ has a section. So without AC projectivity is not
equivalent to being a direct summand of a free module.} (one
cannot even prove the existence of a $\BQ$-linear
embedding of $\BR$ into a free $\BQ$-module $F$, for given 
such an embedding and using a $\BQ$-linear retraction
$F\to\BQ$ one would get a splitting $s:\BR/\BQ\to\BR$ of the
exact sequence $0\to\BQ\to\BR\to\BR/\BQ\to 0$ and therefore
a non-measurable subset $s(\BR/\BQ )\subset\BR$, but it is
known \cite{So} that the existence of such a subset cannot be
proved in set theory without AC).

(ii) Instead of property (c) from Theorem \ref{Mit.2}
the authors of \cite{RG} used a slightly different one,
which is harder to formulate. Probably their property has some
technical advantages.

(iii) Here are some more comments regarding the work
\cite{RG}. First, there is no evidence that the authors of
\cite{RG} knew that Theorem \ref{Mit.1} had been
conjectured by Grothendieck. Second, their notion of
Mittag-Leffler module and their results on infinitely
generated projective modules were probably largely forgotten
(even though they deserve being mentioned in algebra
textbooks). Probably they were ``lost" among many other
powerful and important results of \cite{RG} (mostly in the
spirit of EGA IV).

\section{Families of Tate vector spaces and the
$K_{-1}$-functor}
\label{Tatesection}

\subsection{A class of topological vector spaces}
\label{Tatesubsection}

We consider topological vector spaces over a discrete field $k$.

\medskip

\noindent {\bf Definition.}
A topological vector space is {\it linearly compact} if it is
the topological dual of a discrete vector space.

\medskip

\noindent {\bf Example:} $k[[t]]\simeq k\times
k\times\ldots= (k\oplus k\oplus\ldots)^*$.

\medskip

A topological vector space $V$ is linearly compact if and
only if it has the following 3 properties:

1) $V$ is complete and Hausdorff,

2) $V$ has a base of neighborhoods of $0$ consisting of 
vector subspaces,

3) each open subspace of $V$ has finite codimension.

\medskip

\noindent {\bf Definition.}
A {\it Tate space} is a topological vector space
isomorphic to $P\oplus Q^*$, where $P$ and $Q$ are discrete.

A topological vector space $T$ is a Tate
space if and only if it has an open linearly compact
subspace.

\medskip

\noindent {\bf Example:} $k((t))$ equipped with its usual
topology (the subspaces $t^nk[[t]]$ form a base of
neighborhoods of 0). This is a Tate space because it is a
direct sum of the linearly compact space $k[[t]]$ and the
discrete space $t^{-1} k[t^{-1}]$, or because
$k[[t]]\subset k((t))$ is an open linearly compact subspace.

\medskip

Tate spaces play an important role in the algebraic geometry of curves
(e.g., the ring of adeles corresponding to an algebraic curve
is a Tate space) and also in the theory of
$\infty$-dimensional Lie algebras and Conformal Field Theory.
In fact, they were introduced by Lefschetz (\cite{L},
p.78--79) under the name of locally linearly compact spaces.
The name ``Tate space'' was introduced by Beilinson because
these spaces are implicit in Tate's remarkable work \cite{T}.
In fact, the approach to residues on curves developed in
\cite{T} can be most naturally interpreted in terms of the
canonical central extension of the endomorphism algebra of a
Tate space, which is also implicit in \cite{T}.

\subsection{What is a family of Tate spaces?}
\label{Tatemoddef} 
Probably this question has not been considered. We suggest
the following answer. In the category of topological modules
over a (not necessarily commutative) ring $R$ we define a
full subcategory of {\it Tate $R$-modules.} If
$R$ is commutative then we suggest to consider Tate
$R$-modules as ``families of Tate spaces''. This viewpoint
is justified by Theorems \ref{2} and \ref{1} below. 

\subsubsection{Definitions}
An {\it elementary Tate $R$-module\,} is a topological 
$R$-module isomorphic to $P\oplus Q^*$, where $P$, $Q$ are
discrete projective $R$-modules ($P$ is a left module, $Q$ is
a right one). A {\it Tate $R$-module\,} is a direct summand of
an elementary Tate $R$-module. A Tate $R$-module $M$ is {\it
quasi-elementary\,} if $M\oplus R^n$ is elementary for some
$n\in\BN$.

\medskip

By definition, a morphism of Tate modules is a
continuous homomorphism. The following lemma is very easy.

\begin{lm}  \label{veasy}
Let $P,Q$ be as in the definition of Tate $R$-module. Then
every morphism $Q^*\to P$ has finitely generated image.
\hfill\qedsymbol
\end{lm}

\subsubsection{Examples} \label{noda2}
1) $R((t))^n$ is an elementary Tate $R$-module.

2) A finitely generated projective $R((t))$-module $M$
has a unique structure of topological
$R((t))$-module such that every $R((t))$-linear morphism
$M\to R((t))$ is continuous. This topology is called the
{\it standard topology\,} of $M$. Clearly $M$ equipped with
its standard topology is a Tate $R$-module. In general, it is
not quasi-elementary. E.g., let $k$ be a field, $R:=\{ f\in k
[x] |f(0)=f(1)\}$ and 
\begin{equation} \label{eqone}
M:=\{u=u(x,t)\in k[x]((t))\,|\, u(1,t)=tu(0,t)\}.
\end{equation}
Then $M$ is a finitely generated 
projective $R((t))$-module which is not quasi-elementary
as a  Tate $R$-module (see \ref{nodal}).

\medskip

\noindent {\bf Remark.} The precise relation between finitely
generated projective $R((t))$-modules and Tate $R$-modules
is explained in Theorem \ref{6.3new} below.

\subsubsection{Lattices and bounded submodules}
A submodule $L$ of a topological $R$-module $M$ is said to be
a {\it lattice\,} if it is open and
$L/U$ is finitely generated for every open submodule
$U\subset L$. A subset of a Tate $R$-module $M$ is {\it
bounded\,} if it is contained in some lattice. A lattice $L$
in a Tate module
$M$ is {\it \coflat\,} if $M/L$ is projective.

\medskip

\noindent {\bf Remarks.} (i)  One can show
that a  lattice $L$ in a Tate module $M$ is \coflat if and
only if $M/L$ is flat.

(ii) In every Tate $R$-module lattices exist and, moreover,
form a base of neighborhoods of $0$. On the other hand, a
Tate $R$-module $M$ has a \coflat lattice if and only if $M$
is elementary.

\begin{theor} \label{Tatecountable}
A Tate $R$-module $M$ has the following properties:

(a) $M$ is complete and Hausdorff;

(b) lattices in $M$ form a base of neighborhoods of $0$;

(c) the functor that associates to a discrete $R$-module
$N$ the group $\Hom (M,N)$ of continuous homomorphisms 
$M\to N$ is exact.

If a topological $R$-module $M$ has a countable base of
neighborhoods of $0$ and satisfies (a)-(c) then it is a
Tate $R$-module.
\end{theor}

Only the last statement of the theorem is nontrivial. The
countability assumption is essential in it.

\medskip 

\noindent {\bf Remark.} If a topological $R$-module $M$
satisfies (b) then (c) is equivalent to the following
property: for every lattice $L$ there is a lattice
$L'\subset L$ such that the morphism $M/L'\to M/L$ admits a
factorization $M/L'\to P\to M/L$ for some projective module
$P$.

\subsubsection{Duality}
The {\it dual\,} of a Tate $R$-module $M$ is defined to be
the right $R$-module $M^*$ of continuous homomorphisms 
$M\to R$ equipped with the topology whose base is formed by
orthogonal complements of open bounded submodules 
$L\subset M$. Then $M^*$ is again a Tate module, and
$M^{**}=M$ (it suffices to check this for elementary Tate
modules).

\subsubsection{Tate modules as families of Tate spaces}
\label{asfam}

\begin{theor} \label{2}
The notion of Tate module over a commutative ring $R$ is
local for the flat topology, i.e., for every faithfully flat
commutative $R$-algebra $R'$ the category of Tate $R$-modules
is canonically equivalent to that of Tate $R'$-modules
equipped with a descent datum.
\end{theor}

The proof is based on the Raynaud--Gruson
technique. 

\begin{theor}  \label{1}
Let $R$ be a commutative ring. Then every Tate $R$-module
$M$ is Nisnevich-locally elementary; in other words, there
exists a Nisnevich covering $\Spec R'\to\Spec R$ such that
$R'\hat\otimes_RM$ has a \coflat lattice $L'$. Moreover, for
every lattice $L\subset M$ one can choose $R'$ and $L'$
so that $L'\supset R'\hat\otimes_RL$.
\end{theor}

The proof is not hard. A close statement (Theorem
\ref{Kvanishing}) will be proved in \ref{Kvanishsubsect}.

\medskip

Let me give the definition of Nisnevich covering. A morphism
$\pi:X\to\Spec R$ is said to be a {\it Nisnevich covering\,}
if it is etale and there exist closed subschemes 
$\Spec R=F_0\supset F_1\supset\ldots\supset F_n=\emptyset$
such that each $F_i$ is defined by finitely many equations
and $\pi$ admits a section over $F_{i-1}\setminus F_{i}$,
$i=1,\ldots,n$.  A morphism $\pi:X\to Y$ is a Nisnevich
covering if for every open affine $U\subset Y$ the morphism
$\pi^{-1}(U)\to U$ is a Nisnevich covering. (If $Y$ is
locally Noetherian then an etale morphism $X\to Y$ is a
Nisnevich covering if and only if it admits a section over
each point of $Y$; this is the usual definition.) The
Nisnevich topology is weaker than etale but stronger than
Zariski. The following table may be helpful:
\begin{center}
\begin{tabular}{|c|c|}   \hline
Topology& Stalks of $\mathcal{O}_X$, $X=$Spec$R$\\ 
\hline 
Zariski& Localizations of $R$  \\ \hline
Nisnevich& Henselizations of $R$   \\ \hline
Etale & Strict henselizations of $R$  \\ \hline
\end{tabular}
\end{center}

\subsubsection{Remarks on Theorem \ref{1}}  \label{noda3}

(i) In Theorem \ref{1} {\it one cannot replace ``Nisnevich''
by ``Zariski''.} E.g., we will see in \ref{nodal} that the
Tate module \eqref{eqone} is not Zariski-locally
elementary. 

(ii) It is easy to show that every quasi-elementary Tate
module over a commutative ring $R$ is Zariski-locally
elementary.

\subsection{Tate $R$-modules and $K_{-1}(R)$}
\label{theclass}
How to see that a Tate $R$-module is not quasi-elementary?
We will assign to each Tate $R$-module $M$ a class 
$[M]\in K_{-1}(R)$  so that $[M]=0$ if and only if $M$ is
quasi-elementary. It is easy to define $[M]$ if one uses the
following definition of $K_{-1}(R)$.
 
\subsubsection{$K_{-1}$ via Calkin category}
\label{Calkincat}
First, introduce the following category $\Call =\Call_R$: its
objects are all $R$-modules and the group
$\Hom_{\Call}(M,M')$ of
$\Call$-morphisms $M\to M'$ is defined by 
\begin{displaymath} 
\Hom_{\Call}(M,M'):=\Hom (M,M')/\Hom_f(M,M'),
\end{displaymath} 
where $\Hom_f(M,M')$ is the group of $R$-linear maps 
$A:M\to M'$ whose image is contained in a finitely  generated
submodule of $M'$. Let $\cC\subset\Call$ be the full
subcategory whose objects are {\it projective\,} modules. 
The idempotent completion\footnote{The idempotent completion
of a category $\cB$ is the category $\cB^{\Kar}$ in which an
object is a pair $(B,p:B\to B)$ with $B\in\cB$ and $p^2=p$,
and a morphism $(B_1,p_1)\to (B_2,p_2)$ is a $\cB$-morphism 
$\varphi :B_1\to B_2$ such that $p_2\varphi p_1=\varphi$.
This construction was explained by P.~Freyd in Exercise B2
of Ch.~2 of \cite{Fr} a few years before Karoubi.} of $\cC$
(a.k.a. the Karoubi envelope of $\cC$)  will be denoted by
$\cC^{\Kar}$ or $\cC^{\Kar}_R$ and will be called the {\it
Calkin category\,} of $R$. Let $\cC_{\aleph_0}\subset\cC$ be
the full subcategory  of countably generated projective
$R$-modules and $\cC^{\Kar}_{\aleph_0}$ its idempotent
completion.

\begin{pro}
Every object of $\cC^{\Kar}$ is stably 
equivalent\footnote{Objects $X,Y$ of an additive category
$\cA$ are said to be {\it stably equivalent\,} if 
$X\oplus Z\simeq Y\oplus Z$ for some $Z\in\cA$.} to an
object of $\cC^{\Kar}_{\aleph_0}$. Two objects of
$\cC^{\Kar}_{\aleph_0}$ are stably equivalent in
$\cC^{\Kar}$ if and only if they are stably equivalent in
$\cC^{\Kar}_{\aleph_0}$. So 
\end{pro}

As a corollary, we see that $K_0(\cC^{\Kar})$
is well-defined\footnote{$K_0$ of an additive category
$\cA$ is defined by the usual universal property. It may
exist even if $\cA$ is not equivalent to a small category,
e.g., $K_0$ of the category of {\it all\,} vector spaces
equals $0$.} (even though $\cC^{\Kar}$ is not equivalent to
a small category), and the morphism 
$K_0(\cC^{\Kar}_{\aleph_0})\to K_0(\cC^{\Kar})$ is an
isomorphism. Now define $K_{-1}(R)$ by 
\begin{equation}    \label{K-1def} 
K_{-1}(R):=K_0(\cC^{\Kar}_R).
\end{equation}

\medskip

\noindent {\bf Remarks.} (i) The above definition of
$K_{-1}$ is slightly nonstandard but equivalent to the
standard ones.

(ii) Define the {\it algebraic Calkin ring} by
\begin{displaymath} 
\Calk (R):=\End_{\cC}R^{(\BN )}:=\End R^{(\BN )}/
\End_f R^{(\BN )}, \quad
R^{(\BN )}:=R\oplus R\oplus\ldots
\end{displaymath} 
($\Calk (R)$ is an algebraic version of the
analysts' Calkin algebra, which is defined to be the
quotient of the ring of continuous endomorphisms of a Banach
space by the ideal of compact operators). If
$P\in\cC^{\Kar}_{\aleph_0}$ then
$\Hom_{\cC^{\Kar}}(P,R^{(\BN )})$ is a finitely generated
projective module over $\Calk (R)$. Thus one gets an
antiequivalence between $\cC^{\Kar}_{\aleph_0}$ and the
category of finitely generated
projective $\Calk (R)$-modules, which induces an isomorphism
\begin{displaymath} 
K_{-1}(R)\iso K_0 (\Calk (R))
\end{displaymath} 

\subsubsection{The class of a Tate $R$-module} \label{Tatecl}
Let $\cT_R$ denote the additive category of Tate $R$-modules.
We will define a functor
\begin{equation}    \label{Fi} 
\Phi :\cT_R\to \cC^{\Kar}_R.
\end{equation}
Let $E_R\subset\cT_R$ be the full subcategory of elementary
Tate modules. One gets a functor $\Psi :E_R\to\cC_R$ by
setting $\Psi (P\oplus Q^*):=P$ (here $P$, $Q$ are discrete
projective modules) and defining
$\Psi (f)\in\Hom_{\cC}(P,P_1)$, 
$f:P\oplus Q^*\to P_1\oplus Q_1^*$, to be the image of the
composition 
$P\mono P\oplus Q^*
\buildrel{f}\over{\longrightarrow}P_1\oplus Q_1^*\epi P_1$
in $\Hom_{\cC}(P,P_1)$ (the equality $\Psi (f'f)=
\Psi (f')\Psi (f)$ follows from Lemma \ref{veasy}). The
functor \eqref{Fi} is defined to be the extension of 
$\Psi :E_R\to\cC_R\subset\cC^{\Kar}_R$ to $\cT_R=E_R^{\Kar}$.

Now define the class $[M]$ of a Tate $R$-module
$M$ by $[M]:=[\Phi (M)]\in K_0(\cC^{\Kar}_R)=K_{-1}(R)$.

\subsubsection{$K_0$ of the category of Tate $R$-modules}
\begin{theor}     \label{kzeroofTate}
(i) A Tate $R$-module has zero class in $K_{-1}(R)$ if and
only if it is quasi-elementary.

(ii) $K_0(\cT_R)$ is well-defined (even though $\cT_R$ is not
equivalent to a small category).

(iii) The morphism $K_0(\cT_R)\to
K_0(\cC^{\Kar}_R)=K_{-1}(R)$ induced by \eqref{Fi} is an
isomorphism.

(iv) Every element of $K_0(\cT_R)=K_{-1}(R)$ can be
represented as the class of $R((t))\otimes_{R[t,t^{-1}]}P$ 
for some finitely generated projective $R[t,t^{-1}]$-module
$P$. 
\end{theor}

\noindent {\bf Remark.} The only nontrivial point of the
proof is the surjectivity of the composition
\begin{equation} \label{the-morphism}
K_0(R[t,t^{-1}])\to K_0(R((t)))\to K_0(\cT_R)\to K_{-1}(R), 
\end{equation}
which is used in the proof of (iii) and (iv) (in fact, to
prove (iii) it suffices to use Theorem \ref{existTate}(a)
below). The surjectivity of \eqref{the-morphism} is a
standard fact\footnote{The surjectivity of
\eqref{the-morphism} is a tautology if one uses the
definition of $K_{-1}$ given by H.~Bass \cite{Ba}. But it is
a theorem if one defines $K_{-1}$ by \eqref{K-1def}.} from
$K$-theory. It is proved by noticing that there is a
canonical section
$K_{-1}(R)\to K_0(R[t,t^{-1}])$, namely  multiplication by
the canonical element of 
$K_1(\BZ [t,t^{-1}])$. 

\subsection{Nisnevich-local vanishing of $K_{-1}$}
\label{Kvanishsubsect}
Theorem \ref{1} is closely related\footnote{More precisely: 
Theorem \ref{Kvanishing} follows from
Theorems \ref{1} and \ref{kzeroofTate}(iii); Theorem
\ref{1} follows from Theorem \ref{Kvanishing} and 
\ref{3}(iii).} to the following theorem, which I was unable to
find in the literature.

\begin{theor}  \label{Kvanishing} 
Let $R$ be a commutative ring. Then every element of
$K_{-1}(R)$ vanishes Nisne\-vich-locally.
\end{theor}

\noindent {\bf Remarks.} (i) According to Example 8.5 of
\cite{We2} (which goes back to L.~Reid's work \cite{Re}), 
it is {\it not true\,} that every element of
$K_i(R)$, $i<-1$, vanishes Nisnevich-locally.

(ii) It is known that $K_{-1}$ commutes with filtering
inductive limits. So Theorem \ref{Kvanishing} is equivalent to
vanishing  of $K_{-1}(R)$ for commutative Henselian rings
$R$. I prefer the above formulation of the theorem because
commutation of $K_{-1}$ with filtering inductive limits is
not immediate if one defines $K_{-1}$ by  \eqref{K-1def},
i.e., via the Calkin category.

\medskip

In the proof of Theorem \ref{Kvanishing} given below we use
the definition of $K_{-1}$ from \ref{Calkincat}, but it is
also easy to prove the theorem using the definition of
$K_{-1}$ given by H.~Bass \cite{Ba}.

\begin{proof}
It suffices to show that if $P$ is an $R$-module\footnote{We
need only the case that $P$ is projective, but projectivity is
not used in what follows.}, $F\subset P$ is a finitely
generated submodule, and $\pi\in\End P$ is such that
$\im (\pi^2-\pi )\subset F$ then after Nisnevich localization
there exists
$\widetilde{\pi}\in\End P$ such that
$\widetilde{\pi}^2=\widetilde{\pi}$ and
$\im (\widetilde{\pi}-\pi )\subset F$. 

The idea is to look at the spectrum of $\pi$. There exists a
monic $f\in R[\lambda]$ such that $f(\pi^2-\pi)$ annihilates
$F$. Then $f(\pi^2-\pi)(\pi^2-\pi)=0$. Put
$g(\lambda):=(\lambda^2-\lambda)f(\lambda^2-\lambda)$, then
there is a unique morphism $R[\lambda]/(g)\to\End P$ such
that $\lambda\mapsto\pi$. Put
$S:=\Spec
R[\lambda]/(g)\subset\Spec R\times\mathbb{A}^1$, then
$S\supset\underline{0}\cup\underline{1}$, where
$\underline{0}=\Spec R\times\{ 0\}$ and
$\underline{1}=\Spec R\times\{ 1\}$.

Suppose we have a decomposition
\begin{equation} \label{decomposition}
S=S_0\amalg S_1,\quad S_i \mbox{ open, } 
S_0\supset\underline{0}, S_1\supset\underline{1}.
\end{equation}
Then we can define
$e\in R[\lambda]/(g)=H^0(S,\cO_S)$ 
by $e|_{S_0}=0$, $e|_{S_1}=1$ and define
$\widetilde{\pi}$ to be the image of $e$ in
$\End P$. 

Claim: {\it a decomposition \eqref{decomposition}
exists {\it Nisnevich-locally} on 
$\Spec R$.} Indeed, according to the table at the end of
\ref{asfam}, it suffices to show that this decomposition
exists if $R$ is Henselian. Let $\bar g\in (R/m)[\lambda]$ be
the reduction of $g$ modulo the maximal ideal $m\subset R$.
To get \ref{decomposition} it suffices to choose a
factorization $\bar g=\bar g_0\bar g_1$ so that 
$\bar g_0,\bar g_1$ are coprime, $\bar g_0(0)=0,\bar
g_1(1)=0$ and then lift it to a factorization $g=g_0g_1$.
\end{proof}

\subsection{The dimension torsor}  \label{dimtorsor}

Let $R$ be commutative. Then it follows from Theorem 8.5 of
\cite{We2} that there is a canonical epimorphism 
$K_{-1}(R)\to
H_\mathrm{et}^1(\operatorname{Spec}R,\mathbb{Z})$, 
so a Tate $R$-module $M$ should define 
$\alpha_M\in
H_\mathrm{et}^1(\operatorname{Spec}R,\mathbb{Z})$. 
We will define $\alpha_M$ explicitly as a class
of a certain $\mathbb{Z}$-torsor $\Dim_M$ on $\Spec R$
canonically associated to $M$. $\Dim_M$ is called ``the
torsor of dimension theories'' or ``dimension torsor''.

\subsubsection{The case that $R$ is a field.} If $M$ is a
Tate vector space over a field $R$ the notion of dimension
torsor is well known.\footnote{I copied the
definition below from \cite{Ka3}, but the notion goes back
at least to the physical concept of ``Dirac sea",  which
many years later became the ``infinite wedge construction"
in the representation theory of infinite-dimensional Lie
algebras.} Notice that if $L\subset M$ is open and linearly
compact  then usually $\dim L=\infty$ and
$\dim(M/L)=\infty$. But for any open linearly compact $L$,
$L'\subset M$ one has the {\it relative dimension}
$d_L^{L'}:=\dim(L'/L'\cap L)-\dim(L/L'\cap L)\in\BZ$.

\medskip

\noindent {\bf Definition.}
A {\it dimension theory} on a Tate vector space $M$ is a
function
$$d:\{\mbox{open linearly compact subspaces }L\subset
M\}\to\BZ$$ 
such that $d(L')-d(L)=d_L^{L'}$.

A dimension theory exists and is unique up to adding
$n\in\mathbb{Z}$. So dimension theories on a Tate space form a
$\mathbb{Z}$-torsor. This is $\Dim_M$.

\medskip

\noindent {\bf Example.}
Let $T$ be a $\BZ$-torsor, let $R^{(T)}$ be the vector space
over a field $R$ freely generated by $T$. Then $\BZ$ acts on
$R^{(T)}$, so $R^{(T)}$ becomes a $R[z,z^{-1}]$-module
(multiplication by $z$ coincides with the action of
$1\in\BZ$). Put
$M:=R((z))\otimes_{R[z,z^{-1}]}R^{(T)}$. Then one has
a canonical isomorphism
\begin{equation} \label{eqtwo}
\Dim_M\iso T:
\end{equation}
to $t\in T$ one associates the dimension theory $d_t$ such
that $d_t (L_t)=0$, where $L_t\subset M$ is the
$R[[z]]$-subspace generated by $t$.

\subsubsection{The general case}

If $M$ is a Tate module and $L\subset L'\subset M$ are \coflat
lattices then $L'/L$ is a finitely generated projective
$R$-module, so if $R$ is commutative then 
$d_L^{L'}:= \rank (L'/L)\in H^0(\Spec R,\BZ )$ is
well-defined.

\medskip

\noindent {\bf Definition.} Let $M$ be a Tate module over a
commutative ring $R$. A {\it dimension theory\,} on $M$ is a 
rule that associates to each $R$-algebra $R'$ and each
\coflat lattice $L\subset R'\hat\otimes_RM$ a locally constant
function $d_L: \Spec R'\to\BZ$ in a way compatible with base
change and so that  $d_{L_2}-d_{L_1}=\rank (L_2/L_1)$ for any
pair of \coflat lattices $L_1\subset L_2\subset
R'\hat\otimes_RM$. Here $R'\hat\otimes_RM$ denotes the
completed tensor product.
 
\medskip

Theorem \ref{1} implies that if the functions $d_L$ with the above
properties are defined for all etale $R$-algebras then there
exists a unique way to extend the definition to all
$R$-algebras. It also shows that dimension theories form
a $\BZ$-torsor for the Nisnevich topology.\footnote{In fact,
the categories of $\BZ$-torsors for the Nisnevich, etale,
fppf, and fpqc topologies are equivalent.} It is called the
{\it dimension torsor\,} and denoted by
$\Dim_M$.

One has a canonical isomorphism 
\begin{equation} \label{additivityofdim}
\Dim_{M_1\oplus M_2}\iso\Dim_{M_1}+\Dim_{M_2}.
\end{equation}
So one gets a morphism 
$K_0(\cT_R)=K_{-1}(R)\to
H_\mathrm{et}^1(\operatorname{Spec}R,\mathbb{Z})$. 
It is surjective. Indeed, let $T$ be a $\BZ$-torsor
on $S:=\Spec R$. Then the free $\cO_S$-module $\cO_S^{(T)}$
generated by the sheaf of sets $T$ is equipped with an action
of $\BZ$, so it is a module over
$\cO_S[z,z^{-1}]$ (multiplication by $z$ coincides with the 
action of $1\in\BZ$). This module is locally free of rank one,
so its global sections form a projective $R[z,z^{-1}]$-module
$R^{(T)}$ of rank 1. Therefore
$R((z))\otimes_{R[z,z^{-1}]}R^{(T)}$ is a Tate $R$-module. Its
dimension torsor is canonically isomorphic to $T$ (cf.
\eqref{eqtwo}). 

\subsubsection{Example} \label{nodal}
Let $M$ be the Tate module \eqref{eqone} over 
$R:=\{ f\in k[x] |f(0)=f(1)\}$. Then the $\BZ$-torsor $\Dim_M$
is nontrivial (its pullback to $S:=\Spec (R\otimes_k\bar k)$
corresponds to the universal covering of $S$). So the class of
$M$ in $K_0(\cT_R)=K_{-1}(R)$ is nontrivial and therefore $M$
is not quasi-elementary. Moreover, it does not become
quasi-elementary after Zariski localization.

\subsubsection{The kernel of the morphism $K_{-1}(R)\to
H_\mathrm{et}^1(\operatorname{Spec}R,\mathbb{Z})$ may be
nonzero} \label{Weibexamples} 
Moreover, this can happen even if $R$ is local. Examples can
be found in \cite{We3}. More precisely, \S6 of \cite{We3}
contains examples of algebras $R$ over a field $k$ such that
$H_\mathrm{et}^1(\operatorname{Spec}R,\mathbb{Z})=0$ but
$K_{-1}(R)\ne 0$. In each of these examples $\Spec R$ is a
normal surface with one singular point $x$. Let $R_x$ denote
the local ring of $x$. According to \cite{We1}, the map
$K_{-1}(R)\to K_{-1}(R_x)$ is an isomorphism, so
$K_{-1}(R_x)\ne 0$.

\subsection{The determinant gerbe}
\label{determinantgerbe}

Given a Tate space $M$ over a field Kapranov \cite{Ka3}
defines  its {\it groupoid of determinant theories\,}.
The definition is based on the notion of relative
determinant of two lattices in a Tate space and goes back to
J.-L.~Brylinski \cite{Br} (and further back to the
Japanese school and \cite{ACK}). If $M$ is a Tate module
over a commutative ring $R$ then rephrasing the definition
from \cite{Ka3} in the obvious way one gets a  sheaf of
groupoids on the Nisnevich topology of $S:=\Spec R$ (details
will be explained in \ref{dim&det}). This sheaf of groupoids
is, in fact, an $\cO_S^{\times}$-gerbe. We call it {\it the
determinant gerbe of $M$.}  Associating the class of this
gerbe to a Tate $R$-module $M$ one gets a morphism
\begin{equation} \label{thegerbe}
K_0(\cT_R)=K_{-1}(R)\to H^2_{\Nis}(S,
\cO_S^{\times}). 
\end{equation}
Probably it is well known to $K$-theorists. One can get the
restriction of \eqref{thegerbe} to 
$\Ker (K_{-1}(R)\to H^1_{\et} (\Spec R,\BZ ))$ (and
possibly the morphism \eqref{thegerbe} itself) from the 
Brown--Gersten--Thomason spectral sequence 
(\cite{TT}, \S10.8). More details on the determinant gerbes
will be given in \S\ref{dim&det}.

\subsection{Co-Sato Grassmannian} \label{twgrassman}
Let $M$ be a Tate module over a commutative ring $R$. The
{\it co-Sato Grassmannian\,} of $M$ is the following functor
$\Gras_M$ from the category of commutative $R$-algebras $R'$
to that of sets: $\Gras_M(R')$ is the set of
\coflat lattices in $R'\hatimes_RM$. Given lattices 
$L\subset M$ and $\tilde L\subset M^*$ let
$\Gras_M^{L,\tilde L}(R')\subset\Gras_M(R')$ be the set of
\coflat lattices in $R'\hatimes_RM$ containing
$R'\hatimes_RL$ and orthogonal to $R'\hatimes_R\tilde L$. The
functor $\Gras_M$ is the inductive limit of the subfunctors
$\Gras_M^{L,\tilde L}$, and these subfunctors form a
filtering family. Theorem \ref{1} easily implies the
following proposition.

\begin{pro}   \label{coSat1}
(i) $\Gras_M^{L,\tilde L}$ is an algebraic space proper and of
finite presentation over $\Spec R$. Locally for the Nisnevich
topology of $\Spec R$ it is a projective scheme over 
$\Spec R$. 

(ii) $\Gras_M$ is an ind-algebraic space ind-proper over
$\Spec R$. 
\end{pro}

\noindent {\bf Remarks.} (a) A standard argument based on the
Pl\"ucker embedding (see \ref{weakcoSat}) shows that if the
determinant gerbe of $M$ is trivial then 
$\Gras_M^{L,\tilde L}$ is projective over
$\Spec R$ and $\Gras_M$ is an ind-projective ind-scheme.

(b) Using Proposition \ref{coSat1} it is easy to prove
ind-representability and ind-properness of {\it the
$\cF$-twisted affine Grassmannian} $\GR_{\cF }$ of a reductive
group scheme $G$ over $R$. Here $\cF$ is a $G$-torsor 
on $\Spec R((t))$ and $\GR_{\cF }$ is the functor that sends a
commutative $R$-algebra $R'$ to the set of extensions of
$\cF\otimes_{R((z))}R'((z))$ to a $G$-torsor over 
$\Spec R'[[z]]$ (up to isomorphisms whose restriction to 
$\cF\otimes_{R((z))}R'((z))$ equals the identity).

\subsection{Finitely generated projective $R((t))$-modules
from the Tate viewpoint.} \label{projandTate}

Theorem \ref{6.3new} below says that a finitely
generated projective $R((t))$-module is the same as a Tate
$R$-module equipped with a topologically nilpotent
automorphism. An endomorphism (in particular, an
automorphism) of a Tate $R$-module $M$ is said to be {\it
topologically nilpotent} if  it satisfies the equivalent
conditions of the next lemma.

\begin{lm} \label{6.1new}
Let $M$ be a Tate $R$-module, $T\in\End M$. Then the
following conditions are equivalent:

(i) $T^n\to 0$ for $n\to 0$ (which means that for
every lattices $L,L'\subset M$ there exists $N$ such that
$T^n L'\subset L$ for all $n>N$);

(ii) there exists a (unique) structure of topological
$R[[t]]$-module on $M$ such that $T$ acts as
multiplication by $t$.
\end{lm}

If $M$ is a finitely generated projective $R((t))$-module 
equipped with its standard topology then multiplication by
$t$ is a topologically nilpotent automorphism of $M$. The
next theorem says that the converse statement is also true.

\begin{theor} \label{6.3new}
Let $M$ be a Tate $R$-module and $T:M\to M$ be a
topologically nilpotent automorphism. Equip $M$ with the
topological $R((t))$-module structure such that $tm=T(m)$
for $m\in M$. Then $M$ is a finitely generated projective
$R((t))$-module, and the topology on $M$ is the standard one.
\end{theor}

\begin{theor} \label{8.4}
Let $R$ be commutative. Then the notion of finitely generated
projective $R((t))$-module is local for the fpqc topology of
$\Spec R$. More precisely, let $R'$ be a faithfully flat
commutative $R$-algebra, $R'':=R'\otimes_RR'$, and let
$f,g:R'((t))\to R''((t))$ be defined by $f(a):=1\otimes a$,
$g(a):=a\otimes 1$; then the category
of finitely generated projective $R((t))$-modules is
canonically equivalent to that of finitely generated
projective $R'((t))$-modules $M'$ equipped with an
isomorphism $R''((t))\otimes_fM'\iso R''((t))\otimes_gM'$
satisfying the usual cocycle condition.
\end{theor}

This is an immediate corollary of Theorems \ref{2} and
\ref{6.3new}. 

\medskip

\noindent {\bf Remark.} If $R$ is of finite type over a field
$k$ and the morphism $\Spec R'\to\Spec R$ is a
Zariski covering then Theorem \ref{8.4} is well known from
the theory of non-archimedian analytic spaces 
\cite{BGR,Be}, which is applicable because $R((t))$ is an
affinoid $k((t))$-algebra in the sense of
\ref{affinoidloops}.

\subsection{The dimension torsor of a projective
$R((t))$-module} 
\label{dimtorsproj}

Let $R$ be a commutative ring. Let $M$ be a finitely
generated projective $R((t))$-module equipped with an
isomorphism $\varphi :\det M\iso R((t))$. If $R$ is a field
then $M$ has an $R[[t]]$-stable lattice; moreover, there
is a lattice $L\subset M$ such that
\begin{equation} \label{latticecond}
R[[t]]L\subset L, \quad \varphi (\det L)=R[[t]].
\end{equation}
So it is easy to see that if $R$ is a field then there is a
unique dimension theory $d_{\varphi}$ on $M$ such that
$d_{\varphi}(L)=0$ for all lattices $L\subset M$ satisfying
\eqref{latticecond}. Therefore if $R$ is any commutative ring
then the $\BZ$-torsor $\Dim_M$ is trivialized over each point
of $\Spec R$.

\begin{pro}    \label{existtrivializ}  
These trivializations come from a (unique) trivialization
$d_{\varphi}$ of the $\BZ$-tor\-sor~$\Dim_M$.
\end{pro}

\medskip

By Proposition \ref{existtrivializ} the morphism
$K_0(R((t)))\to H^1_{\et} (\Spec R,\BZ )$ that sends the
class of a projective $R((t))$-module $M$ to the class of
$\Dim_M$ annihilates the kernel of the epimorphism 
$\det :K_0(R((t)))\epi\Pic R((t))$, so we get a morphism
\begin{equation} \label{f}
f:\Pic R((t))\to H^1_{\et}(\Spec R,\BZ )
\end{equation}
such that the diagram
\begin{equation} \label{diagram}
         \begin{array}{ccc}
K_0(R((t)))&\buildrel{\det}\over{\longrightarrow}&
\Pic R((t))\\
\downarrow&\swarrow\\ 
H^1_{\et}(\Spec R,\BZ )\\
     \end{array}
\end{equation}
commutes. The composition
\begin{equation} \label{g}
g:\Pic R[t,t^{-1}])\to
\Pic R((t))\buildrel{f}\over{\longrightarrow} 
H^1_{\et}(\Spec R,\BZ )
\end{equation}
was studied in \cite{We2}. 

\medskip

\noindent {\bf Remarks.} (i)  As explained in \cite{We2}, the
kernels of \eqref{f} and \eqref{g} 
may be
nontrivial (even if $R$ is Henselian). Example: if $k$ is a
field and $R$ is either $k[x^2,x^3]\subset k[x]$ or the
Henselization of $k[x^2,x^3]$ at the singular point of its
spectrum then $\Ker f\simeq k((t))/k[[t]]$, 
$\Ker g\simeq k[t,t^{-1}]/k[[t]]$ (e.g., to show that 
$\Ker f\simeq k((t))/k[[t]]$ for $R=k[x^2,x^3]$ notice that a
line bundle on $\Spec R((t))$ is the same as a triple
consisting of a line bundle on $\Spec k[x]((t))$, a line
bundle on $\Spec k((t))$ and an isomorphism between their
pullbacks to $\Spec k[x]((t))/(x^2)$). It is also explained
in \cite{We2} that $g$ has a splitting (and therefore $f$
has).  Indeed, $\Pic R[t,t^{-1}]=H^1_{\et}(\Spec R, C)$,
where $C$ is the derived direct image of the etale sheaf of
invertible functions on $\Spec R[t,t^{-1}]$, and the
morphism $\BZ\to C$ defined by $n\mapsto t^n$ gives a
splitting. 

(ii) The interested reader can easily lift the diagram
\eqref{diagram} of abelian groups to a commutative diagram of
appropriate Picard groupoids (in the sense of \S 1.4 of
\cite{Del73}).

\section{Almost projective and \ALMOST
projective modules.}  \label{Almostsect}

\subsection{Main definitions and results} \label{APsect}
Recall that every Tate $R$-module has a lattice but not
necessarily a \coflat one. If $M$ is a Tate $R$-module and
$L\subset M$ is a lattice (resp. a bounded open
submodule) then $M/L$ is \ALMOST projective (resp. \allmost
projective) in the sense of the following definitions.

\medskip

\noindent {\bf Definitions.}
An 
{\it elementary \allmost projective $R$-module\,} is a
module  isomorphic to a direct sum of a projective
$R$-module and a finitely generated one. An 
{\it \allmost projective $R$-module\,} is a
direct summand of an elementary \allmost projective module. An
\allmost projective $R$-module $M$ is {\it
quasi-elementary\,} if $M\oplus R^n$ is elementary for some
$n\in\BN$. 

\medskip

\noindent {\bf Definition.} An $R$-module $M$ is \ALMOST 
projective if it can be represented as a direct summand of
$P\oplus F$ with $P$ a projective $R$-module and $F$ an
$R$-module of finite presentation.

\medskip
In fact, there is a reasonable notion of $n$-\allmost
projectivity for any positive $n$, see Remark~3 at the end
of this subsection.

\medskip

\noindent {\bf Remark.} 
It is easy to show that an \allmost projective
module $M$ is quasi-elementary 
if and only if it can be represented as $P/N$ with $P$
projective and $N\subset P$ a
submodule of a finitely generated submodule of $P$.
It is also easy to show that for $P$ and $N$ as above $P/N$ is
\ALMOST projective if and only if $N$ is finitely generated.

\begin{theor} \label{existTate}
(a) Every \allmost projective $R$-module $M_0$ can be
represented as $M/L$ with $M$ being a Tate $R$-module and
$L\subset M$ a bounded open submodule. 

(b) If $M_0$ is \ALMOST projective then in such a
representation $L$ is a lattice.
\end{theor}

\begin{theor}   \label{3}
(i) The notion of \allmost projective module over a
commutative ring $R$ is local for the flat topology, i.e.,
for every faithfully flat commutative $R$-algebra $R'$ 
\allmost projectivity of an $R$-module $M$ is equivalent to
\allmost projectivity of the $R'$-module $R'\otimes_RM$. The
same is true for \ALMOST projectivity.

(ii) For every \allmost projective module $M$ over a
commutative ring $R$ there exists a Nisnevich covering 
$\Spec R'\to\Spec R$ such that $R'\otimes_RM$ is elementary.

(iii) For every quasi-elementary \allmost projective
module $M$ over a commutative ring $R$ there exists a
Zariski covering $\Spec R=\bigcup_i\Spec R_{f_i}$ such that
$R_{f_i}\otimes_RM$ is elementary for all $i$.
\end{theor}

The proof of (i)
is based on the Raynaud--Gruson technique. The proofs of (ii)
and (iii) are much easier. In particular, (iii) easily
follows from Kaplansky's Theorem \cite{Ka}, which says that a
projective module over a local field is free (even if it is
not finitely generated!).

\medskip

\noindent {\bf Remarks.} 1) In statement (ii) of the
theorem {\it one cannot replace ``Nisnevich'' by
``Zariski''.} E.g., the quotient of the Tate $R$-module
\eqref{eqone} by any open bounded submodule is an \allmost
projective module which is not Zariski-locally elementary
(because the Tate module \eqref{eqone} is not, see
\ref{nodal}).

2) My impression is that statement 
(ii) is more important than (i) even though it is much
easier to prove. Statement (i) gives you a peace of mind
(without it one would have two candidates for the notion of
\allmost projectivity), but in the examples of \allmost
projective modules that I know one can prove \allmost
projectivity directly rather than showing that the property
holds locally. The roles of Theorems \ref{2} and \ref{1} in
the theory of Tate $R$-modules are similar. 

3) Although we do not need it in the rest of this work, let
us define the notion of $n$-\allmost projectivity for any
$n\in\BN$: an $R$-module $M$ is {\it $n$-\allmost
projective\,} if in the derived category of $R$-modules $M$
can be represented as a direct summand of $P\oplus F^{\bcdot}$
with $P$ being a projective $R$-module and $F^{\bcdot}$ being
a complex of projective $R$-modules such that $F^i=0$ for
$i>0$ and $F^i$ is finitely generated for
$i>-n$.\footnote{One can show that this definition is
equivalent to the following one: an $R$-module $M$ is {\it
$n$-\allmost projective\,} if a projective resolution of $M$
viewed as a complex in the Calkin category $\cC^{\Kar}_R$ from
\ref{Calkincat} is homotopy equivalent to a direct sum of an
object of $\cC^{\Kar}_R$ and a complex $C^{\bcdot}$ in
$\cC^{\Kar}_R$ such that $C^i\ne 0$ only for $i\le -n$.} One can show that for
$n=1,2$ this is equivalent to the above definitions of
\allmost projectivity and 2-\allmost projectivity and that if
$n>2$ then an
$R$-module $M$ is
$n$-\allmost projective if and only if it is 2-\allmost
projective and for some (or for any) epimorphism $f:P\epi M$
with $P$ projective $\Ker f$ is $(n-1)$-\allmost projective.
One can also show that a module $M$ over a commutative ring is
$n$-\allmost projective if and only if it can
be Nisnevich-locally represented as a direct sum of a
projective module and a module $M'$ having a resolution
$P_{n-1}\to P_{n-2}\to\ldots P_0\to M'\to 0$ by finitely
generated projective modules.

\subsection{Class of an \allmost projective module in
$K_{-1}$}
\label{apclass}
In \ref{Calkincat} we defined the category $\Call$ and its
full subcategory $\cC$ formed by projective modules. Let
$\cC^{\WP}\subset\Call$ denote the full subcategory of
\allmost projective modules. By definition, an \allmost
projective module $M$ is a direct summand of 
$F\oplus P$ with $F$ finitely generated and $P$ projective,
so $M$ viewed as an object of $\cC^{\WP}$ becomes a direct
summand of $P\in\cC$. So we get a fully faithful functor
$\Phi :\cC^{\WP}\to\cC^{\Kar}$ (in fact, it is not hard to
prove that $\Phi$ is an equivalence). To an \allmost
projective
$R$-module $M$ one associates an element $[M]\in
K_{-1}(R):=K_0(\cC^{\Kar})$, namely $[M]$ is the class of
$\Phi (M)\in \cC^{\Kar}$. 

Let $T$ be a Tate $R$-module and $L\subset T$ an open
bounded submodule (so $T/L$ is \allmost projective).
Then $[T/L]=[T]$.

\subsection{The dimension torsor of an \allmost
projective module}   \label{APdimtors}

To an \allmost projective module one associates its
{\it dimension torsor}. The definition is given below.
It is parallel to the definition of the  dimension torsor
of a Tate $R$-module, but there is one new feature:
the dimension torsor of an \allmost projective module is
equipped with a {\it canonical upper semicontinuous
section.}

A submodule $L$ of an \allmost projective $R$-module
$M$ is said to be a {\it lattice\,} if it is finitely
generated. In this case $M/L$ is also \allmost projective. A
lattice $L\subset M$ is said to be {\it \coflat\,} if $M/L$ is
projective\footnote{One can show that if $M$ is \ALMOST
projective this is equivalent to $M/L$ being flat.}. One
shows that  in this case $M/L$ is projective and
$L$ has finite presentation, so \coflat lattices exist if and
only if $M$ is elementary. 

Now let $R$ be commutative. We define a {\it dimension
theory\,} (resp. {\it upper semicontinuous dimension
theory\,}) on an \allmost projective $R$-module $M$ to
be a  rule that associates to each $R$-algebra $R'$ and
each \coflat lattice $L\subset R'\otimes_RM$ a locally
constant (resp. an upper semicontinuous) function 
$d_L: \Spec R'\to\BZ$ in a way compatible with base change
and so that $d_{L_2}-d_{L_1}=\rank (L_2/L_1)$ for any pair of
\coflat lattices $L_1\subset L_2\subset R'\otimes_RM$. The
notion of dimension theory (or upper semicontinuous dimension
theory) does not change if one considers only etale
$R$-algebras instead of arbitrary ones. Dimension theories
on an \allmost projective $R$-module $M$ form a
$\BZ$-torsor for the Nisnevich topology of $\Spec R$, which
is denoted by $\Dim_M$. One defines the {\it canonical upper
semicontinuous dimension theory\,} $d^{\can}$ on $M$ by 
$d^{\can}_ L(x):=\dim_{K_x}(K_x\otimes_{R'}L)$, where $R'$ is
an $R$-algebra, $L\subset R'\otimes_RM$ is a \coflat
lattice, $x\in\Spec R'$, and $K_x$ is the residue field of
$x$. An upper semicontinuous dimension theory on $M$ is the
same as an {\it upper semicontinuous section\,} of $\Dim_M$,
by which we mean a $\BZ$-antiequivariant morphism from the
$\BZ$-torsor $\Dim_M$ to the sheaf of upper semicontinuous
$\BZ$-valued functions on $\Spec R$. Clearly $d^{\can}$ is a
true (i.e., locally constant) section of $\Dim_M$ if and only
if the quotient of $M$ modulo the nilradical $I\subset R$ is 
projective over $R/I$. In this case $d^{\can}$ defines a 
trivialization of $\Dim_M$.

If $N$ is a Tate $R$-module and $L\subset N$ is an open
bounded submodule then the dimension torsor of the \allmost
projective module $N/L$ canonically identifies with that of
$N$.

\section{Finer points: determinants and the canonical
central extension}
\label{dim&det}

Subsection \ref{centralproj} (in which we discuss the
canonical central extension of the automorphism group of an
\allmost projective module) is the only part of this section
used in the rest of the article, namely in
\S\ref{polydisksection}. Therefore some readers
(especially those interested primarily in spaces of
formal loops and refined motivic integration) may prefer to
skip this section. But it contains an interesting (though
slightly vague) picture, which I learned from
A.~Beilinson (see \ref{Beilvision}).

In \ref{words}-\ref{Fermi} we follow \S2 of \cite{BBE}.
In particular, we combine the dimension torsor
and the determinant gerbe into a single object, which is a
Torsor over a certain Picard groupoid (these notions are
defined below). The reason why it is convenient and maybe
necessary to do this is explained in \ref{zachem}. Our terminology is slightly different
from that of \cite{BBE}, and our determinant Torsor is inverse
to that of \cite{BBE}.

\subsection{Terminology}  
\label{words}
According to \S 1.4 of \cite{Del73}, a {\it Picard
groupoid\,} is a symmetric monoidal category $\cA$ such that
all the morphisms of $\cA$ are invertible and the semigroup
of isomorphism classes of the objects of $\cA$ is a group. A
Picard groupoid is said to be {\it strictly commutative\,}
if for every $a\in\Ob\,\cA$ the commutativity isomorphism 
$a\otimes a\iso a\otimes a$ equals $\id_a$. As explained in
\S 1.4 of \cite{Del73}, there is also a notion of  {\it
sheaf of Picard groupoids\,}  (champ de cat\'egories de
Picard) on a site. 

We will work with the following simple examples.

\medskip

\noindent {\bf Examples.} For a commutative ring $R$ we have
the Picard groupoid $\cP ic_R$ of invertible $R$-modules and
the Picard groupoid $\cP ic^{\BZ}_R$ of $\BZ$-graded
invertible $R$-modules (the latter is not strictly commutative
because we use the ``super'' commutativity constraint 
$a\otimes b\mapsto (-1)^{p(a)p(b)}b\otimes a$). For
a scheme $S$ denote by $\cP ic^{\BZ}_S$ 
(resp. $\cP ic_S$) the sheaf of Picard groupoids on the
Nisnevich site of $S$ formed by
$\BZ$-graded invertible $\cO_S$-modules (resp. plain
invertible $\cO_S$-modules, a k.a. $\cO_S^{\times}$-torsors).

\medskip

We need more terminology. An {\it Action\,} of a monoidal
category $\cA$ on a category $\cC$ is a monoidal functor from
$\cA$ to the monoidal category $\Funct (\cC,\cC )$ of
functors $\cC\to\cC$. Suppose $\cA$ acts on $\cC$ and $\cC'$,
i.e., one has monoidal functors 
$\Phi :\cA\to\Funct (\cC,\cC )$ and 
$\Phi' :\cA\to\Funct (\cC',\cC' )$. Then an 
{\it $\cA$-functor\, } $\cC\to\cC'$  is a functor 
$F:\cC\to\cC'$ equipped with isomorphisms 
$F\Phi (a)\iso\Phi'(a)F$ satisfying the natural
compatibility condition (the two ways of constructing an
isomorphism
$F\Phi (a_1\otimes a_2)\iso\Phi'(a_1\otimes a_2)F$
must give the same result). An {\it $\cA$-equivalence\, }
$\cC\to\cC'$ is an $\cA$-functor $\cC\to\cC'$ which is an
equivalence. 

There is also an obvious notion of Action of a sheaf of
monoidal categories $\cA$ on a sheaf of categories $\cC$,
and given an Action of $\cA$ on $\cC$ and $\cC'$ there is
an obvious notion of $\cA$-functor $\cC\to\cC'$ and
$\cA$-equivalence $\cC\to\cC'$.

\medskip

\noindent {\bf Definition.} Let $\cA$ be a sheaf of Picard
groupoids on a site. A sheaf of categories $\cC$ equipped
with an Action of $\cA$ is an {\it $\cA$-Torsor\,} if it is
locally $\cA$-equivalent to $\cA$.
 
\medskip

\noindent {\bf Remark.} The notion of Torsor makes sense even
if $\cA$ is non-symmetric. But $\cA$ has to be symmetric if
we want to have a notion of product of $\cA$-Torsors.

\subsection{The determinant Torsor}  
\label{determinanttorsor} 
Let $R$ be a commutative ring,
$S:=\Spec R$. Slightly modifying the construction of
\cite{Ka3}, we will associate a Torsor over $\cP ic^{\BZ}_S$
to an \allmost projective $R$-module $M$. Recall that a
\coflat lattice $L\subset M$ is a finitely generated
submodule such that $M/L$ is projective. The set of \coflat
lattices $L\subset M$ will be denoted by $G(M)$. In general,
$G(M)$ may be empty, and it is not clear if every
$L_1,L_2\in G(M)$ are contained in some
$L\in G(M)$. But it follows from Theorem \ref{3}(ii)   that
these properties hold after Nisnevich localization (to show
that every $L_1,L_2\in G(M)$ are Nisnevich-locally contained
in some \coflat lattice apply statement (ii) or (iii)
of Theorem \ref{3} to $M/(L_1+L_2)$). In other
words, for every $x\in\Spec R$ the inductive limit of
$G(R'\otimes_RM)$ over the filtering category of all etale
$R$-algebras $R'$ equipped with an $R$-morphism $x\to\Spec
R'$  is a non-empty directed set. 

For each pair $L_1\subset L_2$ in
$G(M)$ one has the invertible $R$-module $\det (L_2/L_1)$. It
is equipped with a $\BZ$-grading (the determinant of an
$n$-dimensional vector space has grading~$n$).

\medskip

\noindent {\bf Definition.} A {\it determinant
theory\,} on $M$ (resp. a {\it weak determinant
theory\,} on $M$)  is a rule $\Delta$ which
associates to each  $R$-algebra $R'$ and each 
$L\in G(R'\otimes_RM)$ an invertible graded $R'$-module 
$\Delta (L)$ (resp. an invertible $R'$-module 
$\Delta (L)$), to each pair
$L_1\subset L_2$ in $G(R'\otimes_RM)$ an isomorphism
\begin{equation} \label{DeltaL1L2}
\Delta_{L_1L_2}: \Delta(L_1)\otimes\det (L_2/L_1)
\iso\Delta(L_2),
\end{equation}
and to each morphism $f:R'\to R''$ of  $R$-algebras
a collection of base change morphisms
$\Delta_f=\Delta_{f,L'}:
\Delta (L' )\to \Delta (R''L' )$, $L'\in G(R'\otimes_RM)$.
These data should satisfy the following conditions:

(i) each 
$\Delta_{f,L'}$ induces an isomorphism 
$R''\otimes_{R'}\Delta (L' )\iso \Delta (R'' L' )$;

(ii) $\Delta_{f_2f_1}=\Delta_{f_2}\Delta_{f_1}$;

(iii)  the
isomorphisms (\ref{DeltaL1L2}) commute with base
change;

(iv) for any triple $L_1\subset L_2\subset L_3$ in
$G(R'\otimes_RM)$ the obvious diagram
\begin{displaymath} 
        \begin{array}{ccc}
\Delta(L_1)\otimes\det(L_2/L_1)\otimes\det(L_3/L_3)&\iso&
\Delta(L_1)\otimes\det(L_3/L_1)\\
\bd&&\bd\\
\Delta(L_2)\otimes\det(L_3/L_2)&\iso&\Delta(L_3)\\
     \end{array}
\end{displaymath} 
commutes.

\medskip

\noindent {\bf Remark.} It follows from Theorem \ref{3}(ii) 
 that the notion of (weak) determinant
theory does not change if one considers only etale
$R$-algebras instead of arbitrary ones.

\medskip

The groupoid of all determinant theories on $M$ is equipped
with an obvious Action of the Picard groupoid $\cP ic^{\BZ}_R$
of invertible $\BZ$-graded $R$-modules: $P\in\cP ic^{\BZ}_R$
sends $\Delta$ to $P\Delta$, where 
$(P\Delta) (L):=P\otimes_R\Delta (L)$.

Determinant theories on $R'\otimes_RM$ for all etale $R$
algebras $R'$ form a sheaf of groupoids $\Det_M$ on the
Nisnevich site of $S:=\Spec R$, which is  equipped with an
Action of the sheaf of Picard groupoids $\cP ic^{\BZ}_S$. It
follows from  Theorem \ref{3}(ii)  that
$\Det_M$ is a Torsor over $\cP ic^{\BZ}_S$. We call it the
{\it determinant Torsor\,} of $M$.

If $M$ is a Tate module (rather than an \allmost projective
one) then the above definition of determinant theory and
determinant Torsor still applies (of course, in this case the
words ``\coflat lattice" should be understood in the sense of
\ref{Tatemoddef} and $\otimes$ should be replaced by
$\hatimes$). If $M$ is an \allmost projective or Tate module
and $L\subset M$ is a lattice then $M/L$ is \allmost
projective and $\Det_{M/L}$ canonically identifies with
$\Det_M$.

\medskip

\noindent {\bf Remark.} Consider the category whose set of
objects is $\BZ$ and whose only morphisms are the identities.
We will denote it simply by $\BZ$. Addition of integers
defines a functor $\BZ\times\BZ\to\BZ$, so $\BZ$ becomes a
Picard groupoid. We have a canonical Picard functor from
$\cP ic^{\BZ}_S$ to the constant sheaf $\BZ$ of Picard
groupoids: an invertible $\cO_S$-module placed in degree $n$
goes to $n$. The $\BZ$-torsor corresponding to the 
$\cP ic^{\BZ}_S$-Torsor $\Det_M$ is the dimension torsor
$\Dim_M$ from \ref{dimtorsor}, \ref{APdimtors}.

\subsection{On the notion of determinant gerbe} 
\label{zachem}

We also have the forgetful functor from the category
of $\BZ$-graded invertible $R$-modules to that of plain
invertible $R$-modules and the corresponding functor
$F:\cP ic^{\BZ}_S\to \cP ic_S$. Notice that $F$ is a
monoidal functor, but {\it not\,} a Picard functor.
Applying $F$ to the $\cP ic^{\BZ}_S$-Torsor
$\Det_M$ one gets a $\cP ic_S$-Torsor, which is the same as
an $\cO_S^{\times}$-gerbe.\footnote{This follows from the
definitions, but also from the Grothendieck-Deligne
dictionary mentioned in \ref{Beilvision} (the complex of
sheaves of abelian groups corresponding via this dictionary
to the sheaf of Picard categories $\cP ic_S$ is
$\cO_S^{\times}[1]$, i.e., $\cO_S^{\times}$ placed in degree
-1).} This is the {\it determinant gerbe\,} considered by
Kapranov
\cite{Ka3} and mentioned in \ref{determinantgerbe}. Its
sections are weak determinant theories. As $F$ does not
commute with the commutativity constraint, there is no
canonical equivalence between the $\cO_S^{\times}$-gerbe
corresponding to a direct sum of \allmost projective modules
$M_i$, $i\in I$, $\Card I<\infty$, and the product of the
$\cO_S^{\times}$-gerbes corresponding to
$M_i$, $i\in I$ (but there is an equivalence which depends
on the choice of an ordering of $I$). This is the source of
the numerous signs in \cite{ACK} and the reason why we
prefer to consider Torsors over $\cP ic^{\BZ}_S$ rather than
pairs consisting of an $\cO_S^{\times}$-gerbe and a
$\BZ$-torsor (as  Kapranov does in
\cite{Ka3}).

\subsection{Fermion modules, determinant theories, and
co-Sato Grassmannian} 
\label{Fermi}
We follow \S\S2.14 -- 2.15 of \cite{BBE}  (in particular, see
Remark (iii) at the end of \S2.15 of \cite{BBE}). 

\subsubsection{Fermion modules and weak determinant theories}
\label{Fermiweak} 
Fix a Tate $R$-module $M$. Let $\Cl (M\oplus M^*)$ denote
the Clifford algebra of $M\oplus M^*$. 
Define a {\it
Clifford module\,} to be a module $V$ over 
$\Cl (M\oplus M^*)$ such that for any $v\in V$ the set
$\{ a\in M\oplus M^*|av=0\}$ is open in $M\oplus M^*$. A
Clifford module $V$ is said to be a {\it fermion
module}\footnote{Motivation of the name: if $M$ is a discrete
projective $R$-module then fermion modules have the form
\linebreak
$(\bigwedge M)\otimes_R\cL$ with $\cL$ being an invertible
$R$-module. Here $\bigwedge M$ is the exterior
algebra of $M$.} if $V$ is fiberwise irreducible  and
projective over $R$. 

If $V$ is a fermion module and $L\subset M$ is a \coflat
lattice let $\Delta_V(L)$ denote the annihilator of 
$L\oplus L^{\perp}\subset M\oplus M^*$ in $V$. As
explained in \cite{BBE}, $\Delta_V(L)$ is a {\it line\,} in
$V$ (i.e.,  a direct summand of $V$ which is an invertible
$R$-module) and $\Delta_V$ is a weak determinant theory:
if $L_1\subset L_2\subset M$ are \coflat lattices then the
isomorphism \eqref{DeltaL1L2} comes from the composition
$\bigwedge^rL_2\to\bigwedge M\to\Cl (M\oplus M^*)$, where $r$
is the rank of $L_2/L_1$ and $\bigwedge M$ is the exterior
algebra of $M$. Thus one gets a functor
$V\mapsto\Delta_V$ from the groupoid of fermion modules to
that of weak determinant theories. As explained in \cite{BBE},
it is an equivalence: to construct the inverse functor
$\Delta\mapsto V_{\Delta}$ one first constructs
$V_{\Delta}$ Nisnevich-locally, then glues the results of
the local constructions, and finally uses Theorem \ref{Mit.1}
to prove that $V_{\Delta}$ is a projective $R$-module.

The equivalences $V\mapsto\Delta_V$ and
$\Delta\mapsto V_{\Delta}$ are compatible with the
Actions of the groupoid $\Pic_R$ of invertible $R$-modules.

\subsubsection{Graded fermion modules and determinant
theories}
As explained in \cite{BBE}, the fermion module $V_{\Delta}$
corresponding to a weak determinant theory $\Delta$ is
equipped with a $T$-grading, where $T$ is the dimension torsor
of $M$. Given a determinant theory on $M$ rather than a weak
determinant theory one gets a $\BZ$-grading on the fermion
module compatible with the $\BZ$-grading of $\Cl (M\oplus
M^*)$ for which $M$ has degree 1 and $M^*$ has degree -1.
Thus $\Det_M$ identifies with the groupoid of $\BZ$-graded
fermion modules.

\subsubsection{The Pl\"ucker embedding of the co-Sato
Grassmannian}  \label{weakcoSat} 
The co-Sato Grassmannian $\Gras_M$ of a Tate $R$-module $M$
was defined in \ref{twgrassman}. Now suppose that the
determinant gerbe of $M$ is trivial and fix a  weak
determinant theory $\Delta$ on $M$. Then we get a line
bundle $\cA_{\Delta}$ on $\Gras_M$ whose fiber over a \coflat
lattice $L$ equals $\Delta (L)$. 

On the other hand, we have the fermion module
$V=V_{\Delta}$ such that $\Delta$=$\Delta_V$ (see
\ref{Fermiweak}). Assigning to a \coflat lattice $L$ the line
$\Delta_V(L)$ one gets a morphism $i:\Gras_M\to\BP$, where
$\BP$ is the ind-scheme of lines in $V$. As explained by
Pl\"ucker, $i$ is a closed embedding.

Clearly $\cA_{\Delta}=i^*\cO (-1)$.

\subsection{A somewhat vague picture} \label{Beilvision}
\subsubsection{The picture I learned from Beilinson}
Let $S$ be a spectrum in the sense of algebraic topo\-logy.
We put $\pi^i(S):=\pi_{-i}(S)$ and define $\tau^{\le k}S$ to
be the  spectrum equipped with a morphism $\tau^{\le k}S\to S$
such that $\pi^i(\tau^{\le k}S)=0$ for $i>k$ and the morphism
$\pi^i(\tau^{\le k}S)\to\pi^i(S)$ is an isomorphism for 
$i\le k$. There is a notion of torsor over a spectrum $S$,
which depends only on $\tau^{\le 1}S$. Namely, an $S$-torsor
is a point of the infinite loop space $L$ corresponding to
$(\tau^{\le 1}S)[1]$ (or equivalently, a morphism from the
spherical spectrum to $S[1]$). A homotopy equivalence between
torsors is a path connecting the corresponding points of $L$,
so  equivalence classes are parametrized by
$\pi^1(S):=\pi_{-1}(S)$.

Beilinson's first remark: an object of the Calkin category
$\cC^{\Kar}_R$ (see \ref{Calkincat}) defines a point of the
infinite loop space corresponding to the $K$-theory spectrum
$K(\cC^{\Kar}_R)$, and as $K(\cC^{\Kar}_R)=K(R))[1]$ it
defines a $K(R)$-torsor. In particular, an \allmost
projective $R$-module $M$ defines a $K(R)$-torsor, whose
class in $\pi_{-1}(K(R))=K_{-1}(R)$ is the class
$[M]$ considered in~\ref{theclass}. If $R$ is commutative then
by Thomason's localization theorem (\cite{TT}, \S10.8)
$K(R)=R\Gamma (S,\cK )$, where $\cK$ is the sheaf of
$K$-theories of $\cO_S$ (this is a sheaf of spectra on the
Nisnevich site of $S$). So the notion of $K(R)$-torsor
should\footnote{Here and in what follows I use the word
``should" to indicate the parts of the picture that I do not
quite understand (probably due to the fact that I have not
learned the theory of sheaves of spectra).} coincide with
that of $\cK$-torsor. Both of them should coincide with that
of $\tau^{\le 1}\cK$-torsor. By Theorem
\ref{Kvanishing}, $\cK^1:=\cK_{-1}=0$, so 
$\tau^{\le 1}\cK=\tau^{\le 0}\cK$ and therefore we get a
morphism
$\tau^{\le 1}\cK=\tau^{\le 0}\cK\to
\cK_{[0,1]}:=\cK^{[-1,0]}:=\tau^{\ge -1}\tau^{\le 0}\cK$. 
So to an \allmost projective $R$-module $M$ there should
correspond a $\cK_{[0,1]}$-torsor $\Delta_M$. According to
Beilinson, {\it $\cK_{[0,1]}$ and $\Delta_M$ should identify
with $\cP ic^{\BZ}_S$ and the Torsor $\Det_M$ from
\ref{determinanttorsor} via the following dictionary,\,} which
goes back to A.~Grothendieck and was used in \S 1.4-1.5 of
\cite{Del73} and in \S4 of \cite{Del87}.

\subsubsection{Grothendieck's Dictionary}
According to it, {\it a Picard groupoid is essentially the
same as a spectrum $X$ with $\pi_i (X)=0$ for $i\ne 0,1$.}
More precisely, the following two constructions become
essentially inverse to each other if the first one is applied
only to infinite loop spaces $X$ with $\pi_i(X)=0$ for $i>1$:

(i) to an infinite loop spaces $X$ one associates its
fundamental groupoid $\Pi (X)$ viewed as a Picard
groupoid\footnote{If $X=\Omega Y$ the group structure on
$\pi_0(X)=\pi_1(Y)$ lifts to a monoidal category structure
on $\Pi (X)$. If $X=\Omega^2Z$ the proof of the
commutativity of $\pi_0(X)=\pi_2(Z)$ ``lifts" to a braiding
on $\Pi (X)$, the ``square of the braiding" map
$t:\pi_0(X)\times\pi_0(X)\to\pi_1(X)$ equals the Whitehead
product $\pi_2(Z)\times\pi_2(Z)\to\pi_3(Z)$, and therefore
$t$ vanishes if $Z$ is a loop space.};

(ii) to a Picard groupoid one associates its classifying
space viewed as an infinite loop space\footnote{The
classifying space $B\cA$ of any symmetric monoidal category
$\cA$ is a $\Gamma$-space (see \cite{Seg}) or if you prefer,
an $E_{\infty}$ space (see \cite{M}). So if every object of
$\cA$ is invertible (i.e., if $\pi_0(B\cA )$ is a group) then 
$B\cA$ is an infinite loop space.}.

For strictly commutative Picard groupoids there is a similar
dictionary and, moreover, a precise reference, namely
Corollary 1.4.17 of \cite{Del73}. The statement from
\cite{Del73} is formulated in a more general context of
sheaves. It says that {\it a sheaf of
strictly commutative Picard groupoids is essentially the same
as a complex of sheaves of abelian groups with cohomology
concentrated in degrees 0 and -1.}

Hopefully, there is also a sheafified version of the
dictionary in the non-strictly commutative case. It should say
that {\it a sheaf $\cA$ of Picard groupoids is essentially the
same as a sheaf of spectra $\cS$ whose sheaves of homotopy
groups $\pi_i$ vanish for $i\ne 0,1$} and that {\it the notion
of
$\cA$-Torsor from \ref{words} is equivalent to that of
$\cS$-torsor.}

\medskip

\subsubsection{}
{\bf Problem:} {\it make the above somewhat vague
picture precise.} The notion of determinant Torsor is very
useful, and its rigorous interpretation in the standard
homotopy-theoretic language of algebraic $K$-theory would be
helpful.

\subsection{The central extension for \allmost projective
modules}   \label{centralproj} 

Let $M$ be an \allmost projective module over a commutative
ring $R$ and $\widetilde{M}$ be the corresponding
quasicoherent sheaf on the Nisnevich topology of $S:=\Spec
R$. Then the sheaf $\AAut M:=\AAut\widetilde{M}$ has a
canonical central extension
\begin{equation} \label{central1}
0\to\cO_S^{\times}\to\hatAut M\to\AAut M\to 0
\end{equation}
Its definition is similar to that of the Tate central
extension of the automorphism group of a Tate vector space
a.k.a. ``Japanese" extension (see \cite{Br,Ka3,PS,BBE}).
Namely, if $M$ has a determinant theory $\Delta$ then 
$\AAut M$ is the sheaf of automorphism of $(M,\Delta )$. This
sheaf does not depend (up to canonical isomorphism) on the
choice of $\Delta$. This allows to define $\hatAut M$ even if
$\Delta$ exists only locally.

Now suppose that a group $R$-scheme $G$ acts on
$M$ (i.e., one has a compatible collection of morphisms 
$G(R')\to\Aut (R'\otimes_RM)$ for all $R$-algebras $R'$).
Then \eqref{central1} induces a canonical central extension
of group schemes
\begin{equation}    \label{central2}
0\to\BG_m\to\widehat{G}\to G\to 0
\end{equation}
(one first defines $\widehat{G}$ as a functor
\{$R$-algebras\}$\to$\{groups\} and then notices that
$\widehat{G}$ is representable because it is a $\BG_m$-torsor
over $G$). If $G$ is abelian we get the commutator map
$G\times G\to\BG_m$. If $M$ is projective (in particular,
if $k$ is a field) the extension \eqref{central2} canonically
splits because in this case there is a canonical determinant
theory on $M$ defined by $\Delta (L)=\det L$. The following
example shows that in general the extension \eqref{central2}
can be nontrivial.

\medskip

\noindent {\bf Example} (A.~Beilinson). Let $k$ be a field, 
$R=k[\varepsilon ]/(\varepsilon^2)$. Fix $g\in k((t))$, then
$L_g:=R[[t]]+\varepsilon gR[[t]]$ is a lattice in the
Tate $R$-module $R((t))$. Put $M:=R((t))/L_g$. Let $G_1$
denote the multiplicative group of $R[[t]]$ viewed as a
group scheme over $R$. On $M$ we have the natural action of
$G_1$ and also the action of $\BG_a$ such that $c\in\BG_a$
acts as multiplication by $1+c\varepsilon g$. So
$G:=G_1\times \BG_a$ acts on $M$. The theory of the Tate
extension (see, e.g., \S3 of \cite{BBE}) tells us that in
the corresponding central extension \eqref{central2} the
commutator of $c\in\BG_a$ and $u\in R[[t]]^{\times}$ equals 
$c\varepsilon\cdot\res (u\cdot dg)$. So the 
extension \eqref{central2} is not  commutative if $dg\ne 0$.

\medskip

\noindent {\bf Remarks.}  (i) To define the central
extension \eqref{central2} it suffices to
have an action of $G$ on $M$ as an object of
the Calkin category $\cC^{\Kar}$ defined in \ref{Calkincat}
(see \S2 of \cite{BBE} for more details). Of course, in this
setting the extension \eqref{central2} may be nontrivial even
if $R$ is a field.

(i\,$'$) One can define the extension \eqref{central2} if
$G$ is any group-valued functor on the category of
$R$-algebras (e.g, a group ind-scheme).

(ii) As explained in \cite{BBE}, the
canonical central extension of the automorphism group of a
Tate vector space should rather be considered as a
``super-extension" (this is necessary to formulate the
compatibility between the extensions corresponding to
Tate spaces $T_1$, $T_2$, and $T_1\oplus T_2$). The same is
true for the canonical central extension of the automorphism
group of an object of the Calkin category $\cC^{\Kar}$. But
in the case of an \allmost projective module $M$ ``super" is
unnecessary because any automorphism of $M$ has degree
0, i.e., preserves the dimension torsor (this follows from
the existence of the canonical upper semicontinuous dimension
theory on $M$, see \ref{APdimtors}).

\section{Applications to spaces of formal loops. ``Refined''
motivic integration}
\label{refinedsection}

In this section all rings and algebras are assumed to be
commutative. We fix a ring $k$. Starting from
\ref{affineloops} we suppose that $k$ is a field.

\subsection{A class of schemes}
We will use the following notation for
affine spaces:
\begin{displaymath} 
\mathbb{A}^I:=\operatorname{Spec} k[x_i]_{i\in I}\, ,\quad
\mathbb{A}^{\infty}:=\mathbb{A}^{\BN}.
\end{displaymath}

We say that a $k$-scheme is {\it nice\,}
if it is isomorphic to $X\times\mathbb{A}^I$, where $X$ is of
finite presentation over $k$ (the set $I$ may be infinite).
An affine scheme is nice if and only if it can be defined by
finitely many equations in a (not necessarily
finite-dimensional) affine space.

\medskip

\noindent {\bf Definition.}
A $k$-scheme $X$ is {\it locally nice\,} (resp. {\it
Zariski-locally nice, etale-locally nice}) if it becomes
nice after Nisnevich localization (resp. Zariski or etale
localization). $X$ is {\it differentially nice\,} if for
every open affine $\Spec R\subset X$ the $R$-module
$\Omega_R^1:=\Omega_{R/k}^1$ is \ALMOST projective.

By Theorem \ref{3}(i) etale-local niceness implies
differential niceness. I do not know if etale-local niceness
implies local niceness. Local niceness does not imply
Zariski-local niceness (see \ref{ex} below).

\medskip

For a differentially nice $k$-scheme $X$ one defines the
{\it dimension torsor\,} $\Dim_X$: if $X$ is an affine
scheme $\Spec R$ then $\Dim_X$ is the dimension torsor of the
\allmost projective $R$-module $\Omega_R^1$, and for a general
$X$ one defines $\Dim_X$ by gluing together the torsors 
$\Dim_U$ for all open affine $U\subset X$. If $X'\subset X$ is
a closed subscheme defined by finitely many equations and $X$
is differentially nice then $X'$ is; in this situation
$\Dim_{X'}$ canonically identifies with the restriction of
$\Dim_X$ to $X'$.

In the next subsection we will see that the dimension
torsor of a locally nice $k$-scheme may be nontrivial, and on
the other hand, there exists a locally nice $k$-scheme  with
trivial dimension torsor which is not Zariski-locally nice.

\subsection{Examples}   \label{ex}
(i) Define 
$i:\mathbb{A}^\infty\to\mathbb{A}^\infty$ by
$i(x_1,x_2,\ldots):=(0,x_1,x_2,\ldots)$. Take
$\mathbb{A}^1\times\mathbb{A}^\infty$ and then glue
$(0,x)\in\mathbb{A}^1\times\mathbb{A}^\infty$ with
$(1,i(x))\in\mathbb{A}^1\times\mathbb{A}^\infty$. Thus one
gets a locally nice $k$-scheme $X$ whose dimension torsor is
nontrivial and even not Zariski-locally trivial.

(ii) Let $M$ be an \allmost projective module
over a finitely generated algebra $R$ over a Noetherian ring
$k$. Let $X$ denote the spectrum of the symmetric algebra of
$M$. Then $X$ is locally nice. This follows from Theorem
\ref{3}(ii) and the next theorem, which is due to H.~Bass 
(Corollary 4.5 from \cite{Ba2}).

\begin{theor} \label{Bass} 
If $R$ is a commutative Noetherian ring whose spectrum is
connected then every infinitely generated projective
$R$-module is free. 
\end{theor} 

It is easy to deduce from Theorem \ref{Bass} that $X$ is
Zariski-locally nice if and only if the class of $M$ in
$K_{-1}(R)$ vanishes locally for the Zariski topology (to
prove the ``only if" statement consider the restriction of
$\Omega^1_X$ to the zero section $\Spec R\mono X$). If
$R$ and $M$ are as in \eqref{eqone} then we get the above
Example (i).

(iii) There exists a locally nice scheme $X$ over a field $k$
which is not Zariski-locally nice but has trivial dimension
torsor.\footnote{In \ref{unpleasanthappen} we will see that
one can get such $X$ from the loop space of a smooth affine
manifold.} According to (ii), to get such an example it
suffices to find a finitely generated $k$-algebra $R$ and an
\allmost projective $R$-module $M$ such that
$H_\mathrm{et}^1(\operatorname{Spec}R,\mathbb{Z})=0$ but the
class of $M$ in $K_{-1}(R)$ is not Zariski-locally trivial.
\S6 of \cite{We3} contains examples of finitely generated
normal $k$-algebras $R$  with $K_{-1}(R)\ne 0$. In each of
them $\Spec R$ has a unique singular point $x$, and according
to \cite{We1}, the map $K_{-1}(R)\to K_{-1}(R_x)$ is an
isomorphism. Now take any nonzero element of $K_{-1}(R)$ and
represent it as a class of an \allmost projective $R$-module
$M$.

\subsection{Generalities on ind-schemes}
\label{classof}

The key notions introduced in this subsection are those of
reasonable, T-smooth, and Tate-smooth ind-scheme (see
\ref{reasonable}-\ref{Tsm}). 

\subsubsection{Definition of ind-scheme and formal scheme.}
\label{glossary}

Functors from the category of $k$-algebras to that of sets
will be called ``spaces". E.g., a $k$-scheme can be
considered as a space. {\it ``Subspace"\,} means
``subfunctor". A subspace $Y\subset X$ is said to be {\it
closed\,} if for every (affine) scheme $Z$ and every 
$f:Z\to X$ the subspace $Z\times_XY\subset Z$ is a closed
subscheme. 

Let us agree that an \emph{ind-scheme} is a space which can
be represented as $\limto X_{\alpha}$, where $\{X_{\alpha}\}
$ is a directed family of quasi-compact schemes such that all
the maps $i_{\alpha\beta}:\,X_{\alpha}\to X_{\beta}$,  
$\alpha\le\beta $, are closed embeddings. (Notice that if the
same space can also be represented as the inductive
limit of a directed family of quasi-compact schemes
$X'_{\beta}$ then each $X'_{\beta}$ is contained in some
$X_{\alpha}$ and each $X_{\alpha}$ is contained in some
$X'_{\beta}$.) If $X$ can be represented as above so that the
set of indices $\alpha$ is countable then $X$ is said to be
an \emph{$\aleph_0$-ind-scheme}. If $P$ is a property
of schemes stable under passage to closed subschemes then we
say that $X $ satisfies the \emph{ind}-$P$ property if each
$X_{\alpha}$ satisfies $P$. E.g., one has the notion of
ind-affine ind-scheme and that of ind-scheme of ind-finite
type. 

Set $X_{\red}:=\limto X_{\alpha\red} $; an ind-scheme $X $ is said to
be \emph{reduced}\/ if $X_{\red}=X  $.

\subsubsection{$\cO$-modules and pro-$\cO$-modules on
ind-schemes} \label{sec:theory.set.p} 

A {\it pro-module\,} over a ring $R$ is defined to be a
pro-object\footnote{A nice exposition of
the theory of pro-objects and ind-objects of a category is
given in \S8 of \cite{GV}. See also the Appendix
of \cite{AM}.} of the category of $R$-modules. 
We identify the category of Tate $R$-modules with a full
subcategory of the category of pro-$R$-modules by associating
to a Tate $R$-module the projective system formed by its
discrete quotient modules. 

An {\it $\cO$-module\,} (resp. a {\it pro-$\cO$-module\,})
$P$ on a space $X$ is a rule that assigns to a commutative
algebra $A$ and a point $\phi\in X(A)$ an $A$-module (resp.
a pro-$A$-module) $P_{\phi}$, and to any morphism of
algebras $f: A \to B$ a $B$-isomorphism
$f_P:\,B\mathop{\otimes}\limits_f P_{\phi}\iso P_{f\phi}$ in
a way compatible with  composition of $f$'s. An $\cO$-module
on a scheme $Y$ is the same as a quasicoherent sheaf of
$\cO_Y$-modules, and an $\cO$-module $P$ on an ind-scheme
$X=\limto X_{\alpha}$ is  the same as a collection of
$\cO$-modules $P_{X_{\alpha}}$ on $X_{\alpha}$
together with identifications
$i_{\alpha\beta}^* P_{X_{\beta}} =P_{X_{\alpha}}$ for
$\alpha\le\beta $  that satisfy the obvious transitivity
property. 

\medskip

A pro-$\cO$-module is said to be a {\it Tate sheaf\,}  if for
each $\phi$ as above the pro-module $P_{\phi}$ is a Tate
module. 

The {\it cotangent sheaf\,} $\Omega^1_X$ of an ind-scheme
$X=\limto X_{\alpha}$ is the pro-$\cO$-module whose
restriction to each $X_{\alpha}$ is defined by the
projective system $i_{\alpha\beta}^*\Omega^1_{X_{\beta}}$,
$\beta\ge\alpha$ (here $i_{\alpha\beta}$ is the embedding
$X_{\alpha}\to X_{\beta}$).

\subsubsection{The notion of reasonable ind-scheme}
\label{reasonable}

The following definitions are due to A.~Beilinson. A closed
quasi-compact  subscheme $Y$ of an ind-scheme  $X$ is called
\emph{reasonable} if for any closed subscheme $Z\subset X$
containing $Y$ the ideal of $Y$ in $\cO_Z$ is finitely
generated. Notice that reasonable subschemes of $X$ form a
directed set. An ind-scheme $X$ is \emph{reasonable} if $X$
is the union of its reasonable subschemes, i.e., if it can
be represented as $\limto X_{\alpha}$, where all
$X_{\alpha}$'s are reasonable.

Any scheme is a reasonable ind-scheme. A closed subspace
of a reasonable ind-scheme is a reasonable ind-scheme.
The product of two reasonable ind-schemes is reasonable. The
completion of any ind-scheme along a reasonable closed
subscheme is a reasonable ind-scheme.

\subsubsection{Main example: ind-scheme of formal loops}
\label{loopspace}
Let $Y$ be an affine scheme over $F:=k((t))$. Define a
functor $\cL Y$ from the category of $k$-algebras to that of
sets by $\cL Y(R):=Y(R\hat\otimes F)$, $R\hat\otimes
F:=R((t))$. It is well known and easy to see that $\cL Y$ is
an ind-affine ind-subscheme. This is the {\it ind-scheme of
formal loops\,} of $Y$. If $Y$ is an affine scheme of {\it
finite type\,} over $F$ then $\cL Y$ is a reasonable
$\aleph_0$-ind-scheme.

\subsubsection{Formal schemes}    \label{formal}
We define a \emph{formal scheme} to be an ind-scheme $X$ such
that $X_{\red}$ is a scheme. An
\emph{$\aleph_0$-formal scheme} is a formal scheme which is
an $\aleph_0$-ind-scheme.  The \emph{completion} of an
ind-scheme $Z$ along a closed subscheme $Y\subset Z$ is the
direct limit of closed subschemes $Y'\subset Z$ such that
$Y'_{\red} =Y_{\red}$. In the case of formal schemes we write
``affine" instead of ``ind-affine". A formal scheme $X$ is
affine if and only if $X_{\red}$ is affine.

\medskip

\noindent {\bf Remark.} As soon as you compare the above
definition of formal scheme with the one from EGA~I you see
that they are not equivalent (even in the affine case) but
the difference is not big: an $\aleph_0$-formal scheme in our
sense which is reasonable  in the sense of \ref{reasonable}
is a formal scheme in the sense of EGA~I, and on the other
hand, a Noetherian formal scheme in the sense of EGA~I is
a formal scheme in our sense.

\subsubsection{Formal smoothness}
Following Grothendieck (\cite{Gr64}, \cite{Gr67}), we say
that $X$ is \emph{formally smooth}\/ if for every
$k$-algebra $A$ and every nilpotent ideal $I\subset A$ the
map $X(A)\to X(A/I)$ is surjective. A morphism
$X\to Y$ is said to be formally smooth if for every
$k$-algebra $k'$ and every morphism $\Spec k'\to Y$ the
$k'$-space $X\times_Y\Spec k'$ is formally smooth. Clearly
formal smoothness of any ind-scheme (resp. a reasonable
ind-scheme) is equivalent to formal smoothness of its
completions along all closed subschemes (resp. all
reasonable closed subschemes).

\begin{theor}    \label{locality}
(i) For reasonable formal schemes formal smoothness is an
etale-local property.

(ii) A reasonable closed subscheme of a formally smooth
ind-scheme is differentially nice.

(iii) If a reasonable $\aleph_0$-ind-scheme $X$ is formally
smooth then $\Omega^1_X$ is
a Tate sheaf (the notions of cotangent sheaf of an
ind-scheme, Tate sheaf, and \MLT sheaf are defined in
\ref{sec:theory.set.p}).
\end{theor}

\medskip
In the case of schemes  statement (i) of the theorem
was proved by Grothendieck (cf. Remark 9.5.8 from \cite{Gr})
modulo the conjecture on the local nature of projectivity
(which was proved a few years later in \cite{RG}). The proof
of Theorem \ref{locality} in the general case  is slightly
more complicated but based on the same ideas.

\subsubsection{T-smoothness and Tate-smoothness}
\label{Tsm} 
We say that a reasonable ind-scheme $X$ is
{\it T-smooth\,} if

(i) every reasonable closed subscheme of $X$ is locally
nice;

(ii) $X$ is formally smooth.

A T-smooth ind-scheme $X$ is said to be {\it Tate-smooth\,} 
if its cotangent sheaf is a Tate sheaf (according to Theorem
\ref{locality}(iii), this is automatic for
$\aleph_0$-ind-schemes).

\medskip

\noindent {\bf Remark.} In the above definitions we do not
require every closed subscheme of $X$ to be
contained in a formally smooth subscheme. It is not
clear if this property holds for $\cL (SL(n))$ or for the
affine Grassmannian, even though these ind-schemes are \IS.
See also Remark (ii) from \ref{affineloops}.

\subsubsection{Dimension torsor}
Let $X$ be a reasonable ind-scheme such that all its
reasonable closed subschemes are differentially nice (by
Theorem \ref{locality}(ii) this is true for any formally
smooth reasonable ind-scheme). Then there is an obvious
notion of the {\it dimension torsor\,} of $X$: for each
reasonable closed subscheme $Y\subset X$ one has the
dimension torsor $\Dim_Y$, and if $Y'\subset Y$ are
reasonable closed subschemes then $\Dim_{Y'}$ identifies
with the restriction of $\Dim_Y$ to $Y'$.

\subsubsection{Relation with the Kapranov-Vasserot theory}
The notion of $T$-smooth ind-scheme is similar to the
notion of ``smooth locally compact ind-scheme" introduced by
M.~Kapranov and E.~Vasserot (see Definition 4.4.4 
from \cite{KV}). Neither of these classes of ind-schemes
contains the other one. The theory of $\cD$-modules on smooth
locally compact ind-schemes developed in \cite{KV} renders
to the class of $T$-smooth ind-schemes, and the same is
true for the Kapranov-Vasserot theory of de Rham complexes
(which goes back to the notion of {\it chiral de Rham
complex\,} from \cite{MSV}). According to A.~Beilinson
(private communication), these theories, in fact, render to
the class of formally smooth reasonable ind-schemes, which
contains both ``smooth locally compact" ind-schemes in the
sense of \cite{KV} and $T$-smooth ones.

\subsection{Loops of an affine manifold}
\label{affineloops}

From now on we assume that $k$ is a field (I have not checked if Theorems
\ref{ind-smooth} and \ref{ind-smooth2} hold for any
commutative ring $k$). So $F=k((t))$ is
also a field. For any affine $F$-scheme $Y$
one has the ind-scheme of formal loops $\cL Y$
(see \ref{loopspace}).

\begin{theor}  \label{ind-smooth}
If an affine $F$-scheme $Y$ is smooth then $\cL Y$ is
\IS.
\end{theor}

\noindent {\bf Remarks.} (i) The theorem is not hard. It is
only property (i) from the definition of T-smoothness (see
\ref{Tsm}) that requires some efforts. See \ref{loopsketch}
for more details. 

(ii) If $Y$ is a smooth affine $F$-scheme then by Theorem
\ref{ind-smooth} every reasonable closed subscheme 
$X\subset \cL Y$ is locally nice. But there exist $Y$ and
$X\subset\cL Y$ as above such that $X$ is not Zariski-locally
nice. One can choose $Y$ and $X$ so that $\Dim_X$ is not
Zariski-locally trivial. But one can also  choose $Y$ and $X$
so that $\Dim_X$ is trivial but $X$ is not Zariski-locally
nice. See \ref{asbadasposs} for examples of these situations.
According to H.~Bass \cite{Ba}, $K_{-1}$ of a regular ring is
zero, so in these examples $\cL Y$ cannot be
represented (even Zariski-locally) as the union of an
increasing sequence of smooth closed subschemes.

\medskip

In the next subsection we formulate an analog of Theorem
\ref{ind-smooth} for affinoid analytic spaces (this is a
natural thing to do in view of \ref{follows}).

\subsection{Loops of an  affinoid space}
\label{affinoidloops}

We will use the terminology from
\cite{BGR} (which goes back to Tate) rather than the one from
\cite{Be}. Let 
$F\langle z_1,\ldots ,z_n\rangle\subset F[[z_1,\ldots z_n]]$
be the algebra of power series which converge in the polydisk
$|z_i|\le 1$. As $F=k((t))$ one has
$F\langle z_1,\ldots ,z_n\rangle =k[z_1,\ldots ,z_n]((t))$.
For every $k$-algebra $R$ the $F$-algebra 
$R\hat\otimes F=R((t))$ is equipped with the norm whose unit
ball is $R[[t]]$. In particular, 
$F\langle z_1,\ldots ,z_n\rangle$ is a Banach algebra. An 
{\it affinoid $F$-algebra\,} is a topological
$F$-algebra isomorphic to a quotient of 
$F\langle z_1,\ldots ,z_n\rangle $ for some $n$. All morphisms
between affinoid $F$-algebras are automatically continuous
(see, e.g., \S6.1.3 of \cite{BGR}). The category of {\it
affinoid analytic spaces\,} is defined to be dual to that of
affinoid $F$-algebras; the affinoid space corresponding to an
affinoid $F$-algebra $A$ will be denoted by $\cM (A)$.

For an affinoid analytic space $Z=\cM (A)$ and a $k$-algebra
$R$ denote by $\cL Z(R)$ the set of continuous
$F$-homomorphisms from $A$ to the Banach $F$-algebra
$R\hat\otimes F=R((t))$. It is easy to see that the functor
$\cL Z$ is a reasonable affine
$\aleph_0$-formal scheme in the sense of \ref{formal},
\ref{reasonable} (and therefore an affine formal scheme in
the sense of EGA~I). E.g., if
$Z$ is the unit disk then
$\cL Z$ is the completion of the ind-scheme of formal Laurent
series along the subscheme of formal Taylor series.

\begin{theor}  \label{ind-smooth2}
If an affinoid space $Z$ is smooth then the formal scheme
$\cL Z$ is \IS. In particular, $(\cL Z)_{\red}$ is a locally
nice scheme. 
\end{theor}

\subsection{Theorem \ref{ind-smooth} follows from Theorem
\ref{ind-smooth2}} \label{follows}
Let $Y=\Spec B$ be a closed subscheme of 
$\BA^n=\Spec F[z_1,\ldots ,z_n]$. The ind-scheme $\cL Y$ is
the union of its closed subschemes $\cL_NY$ defined by
$(\cL_NY)(R):=Y(R)\cap (t^{-N}R[[t]])^n\subset R((t))^n$
for any $k$-algebra $R$. The completion of $\cL Y$ along
$\cL_NY$ identifies with $\cL Y_N$, where $Y_N$ is the
affinoid analytic space defined by
\begin{displaymath} 
Y_N:=\cM (B_N),\quad B_N:=B\otimes_{F[z_1,\ldots ,z_n]}
F\langle t^Nz_1,\ldots ,t^Nz_n\rangle
\end{displaymath}
(in other words, $Y_N$ is the intersection of $Y$ with the
polydisk of radius $r^n$, $r:=|t^{-1}|>1$). Therefore
Theorem
\ref{ind-smooth} follows from Theorem
\ref{ind-smooth2}.

\subsection{Sketch of the proof of Theorem \ref{ind-smooth2}}
\label{loopsketch}
The formal smoothness of $\cL Z$ immediately follows from the
definitions. It is also easy to describe the cotangent sheaf
of $\cL Z$. Let $A$ be the affinoid
$F$-algebra corresponding to $Z$. Every
finite-dimensional vector bundle $E$ on
$Z$ defines a Tate sheaf $\cL E$ on $\cL Z$: if 
$\Spec R\subset\cL Z$ is a closed affine subscheme and
$f:A\to R((t))$ corresponds to the morphism 
$\Spec R\mono\cL Z$ then the pullback of $\cL E$ to
$\Spec R$ is the Tate $R$-module 
$R((t))\otimes_A\Gamma (Z,E)$. The proof of the next lemma
is straightforward.

\begin{lm}    \label{easy}
The cotangent sheaf of $\cL Z$ identifies with the Tate
sheaf $\cL\Omega^1_Z$ corresponding to the cotangent bundle
$\Omega^1_Z$ of the analytic space $Z$.
\end{lm}

\noindent {\bf Corollary.} {\it Let $\Spec R\subset\cL Z$ be a
reasonable closed subscheme. Let $M$ be the module of global 
sections of the pullback to $\Spec R$ of the Tate
sheaf $\cL\Omega^1_Z$. Then $\Omega^1_R$ is the quotient of
the Tate $R$-module $M$ by some lattice.}

\medskip

It remains to show that a reasonable  closed
subscheme $\Spec R\subset\cL Z$ is locally nice. It easily
follows from the above Corollary and Theorem \ref{1} that
after Nisnevich localization $\Omega_R^1$ becomes a
direct sum of a free module and a module of finite
presentation. This is a linearized version of local
niceness. To deduce local niceness from its linearized version
one works with the implicit function theorem.

\subsection{The renormalized dualizing complex}
\label{renorm}

Fix a prime $l\ne\Char k$. Let $D^b_c(X,\BZ_l)$ denote the
appropriately defined bounded constructible $l$-adic derived
category on a scheme $X$ (see \cite{E,Ja}). For a general
locally nice $k$-scheme $X$ there is no natural way to define
the dualizing complex $K_X\in D^b_c(X,\BZ_l)$. Indeed, if
$X$ is the product of $\BA^{\infty}$ and a $k$-scheme $Y$ of
finite type and if $\pi :X\to Y$ is the projection then $K_X$
should equal
$\pi^*K_Y\otimes(\BZ_l[2](1))^{\otimes\infty}$, which
makes no sense. But suppose that the dimension $\BZ$-torsor
$\Dim_X$ is trivial and that we have chosen its
trivialization $\eta$. Then one can define the {\it
renormalized dualizing complex\,} 
$K_X^{\eta}\in D^b_c(X,\BZ_l)$. The definition (which is
straightforward) is given below. The reader can skip it and
go directly to~\ref{RGammac}.

First assume that $X$ is nice, i.e., there exists a morphism
$\pi: X\to Y$ such that $Y$ is a $k$-scheme of finite type
and $X$ is $Y$-isomorphic to $Y\times\BA^I$ for some set $I$. 
Let $\cC_X$ be the category of all such pairs $(Y,\pi )$. A
morphism $f:(Y,\pi )\to (Y',\pi')$ is defined to be a
morphism $f:Y\to Y'$ such that $\pi'=f\pi$. Such $f$ is
unique if it exists. The category $\cC_X$ is equivalent to a
directed set. So to define $K_X^{\eta}$ it suffices to define
a functor 
\begin{equation} \label{thefunctor}
\cC_X\to D^b_c(X,\BZ_l), \quad (Y,\pi )\mapsto K_X^{\eta
,\pi} 
\end{equation}
which sends all morphisms to isomorphisms.

If $(Y,\pi )\in\cC_X$ then $\pi^*\Omega^1_Y\subset\Omega^1_X$
is locally of finite presentation and
$\Omega^1_X/\pi^*\Omega^1_Y$ is locally free. So for every
open affine $U\subset X$ one has the \coflat lattice
$\Gamma (U,\pi^*\Omega^1_Y)\subset\Gamma (U,\Omega^1_X)$ and
therefore a section of the torsor $\Dim_X$ over $U$. These
sections agree with each other, so we get a global section
$\eta_{\pi}$ of $\Dim_X$. Put
\begin{equation} \label{m}
m:=\eta_{\pi}-\eta\in H^0(X,\BZ ), 
\end{equation}
\begin{equation}
K_X^{\eta ,\pi}:=\pi^*K_Y\otimes (\BZ_l[2](1))^{\otimes m}, 
\end{equation}
Now let $f:(Y,\pi )\to (Y',\pi')$ be a morphism. One easily
shows that $f:Y\to Y'$ is smooth\footnote{Choosing a section
$Y\to X$ one sees that $Y$ is $Y'$-isomorphic to a retract of
$Y'\times\BA^J$ for some $J$. So $f$ is formally smooth and
therefore smooth.}, so one has a canonical isomorphism
\begin{equation} \label{standard}
K_Y\iso f^*K_{Y'}\otimes (\BZ_l[2](1))^{\otimes d},
\end{equation}
where $d$ is the relative dimension of $Y$ over $Y'$. It is
easy to see that $\pi^*d=\eta_{\pi'}-\eta_{\pi}$, so
\eqref{standard} induces an isomorphism 
$\alpha_f:K_X^{\eta ,\pi}\to K_X^{\eta ,\pi'}$. We define
\eqref{thefunctor} on morphisms by $f\mapsto\alpha_f$.

So we have defined $K_X^{\eta}$ if $X$ is nice. The formation
of $K_X^{\eta}$ commutes with etale localization of $X$. It is
easy to see that $\Ext^i (K_X^{\eta}, K_X^{\eta})=0$ for
$i<0$. So by Theorem 3.2.4 of \cite{BBD} there is a unique way
to extend the definition of $K_X^{\eta}$ to all etale-locally
nice $k$-schemes $X$ so that the formation of $K_X^{\eta}$
still commutes with etale localization.

\subsection{$R\Gamma_c$ of a locally nice scheme}
\label{RGammac}

Suppose we are in the situation of \ref{renorm}, i.e., we
have a locally nice $k$-scheme $X$, a trivialization $\eta$ of
its dimension torsor, and a prime $l\ne\Char k$. Assume
that $X$ is quasicompact and quasiseparated. Then we put 
\begin{equation} \label{RGammaceta}
R\Gamma_c^{\eta}(X\otimes\bar k,\mathbb{Z}_l):=
R\Gamma(X\otimes_k\bar k,K_X^{\eta})^*,
\end{equation}
where $K_X^{\eta}$ is the renormalized dualizing complex
defined in \ref{renorm}.  
$R\Gamma_c^{\eta}(X\otimes\bar k,\mathbb{Z}_l)$ is an object
of $D^b_c(\Spec k,\BZ_l)$, i.e., of the appropriately defined
bounded constructible derived category of $l$-adic
representations of $\Gal (k^s/k)$, where $k^s$ is a separable
closure of~$k$.

\medskip

\noindent {\bf Problems.} 1) Define an object of the
triangulated category of $k$-motives \cite{VSF,VMW} whose
$l$-adic realization equals
$R\Gamma_c^{\eta}(X\otimes\bar k,\mathbb{Z}_l)$ for each
$l\ne\Char k$ (Voevodsky says this can probably be done).

2) Now suppose that the determinant gerbe of $X$ is
trivial and we have fixed its trivialization $\xi$. 
Can one canonically lift 
$R\Gamma_c^{\eta}(X\otimes\bar k,\mathbb{Z}_l)$ to an object
of the motivic stable homotopy category depending on $\eta$ 
and $\xi$? Or at least, can one canonically lift 
$R\Gamma_c^{\eta}(X\otimes\bar k,\mathbb{Q}_l)$ to an object
of the motivic stable homotopy category tensored by $\BQ$?
(The motivic stable homotopy category a.k.a.
$\BA^1$ stable homotopy category was defined in
\cite{Vo}). Reason why $\xi$ is supposed to exist and to be
fixed: if $k=\BR$ this allows to define 
$R\Gamma_c^{\eta ,\xi}(X(\BR ),\mathbb{Z})$.

\subsection{``Refined'' motivic integration}
\label{refined}

Suppose that in the situation of Theorem \ref{ind-smooth2}
the canonical bundle $\det\Omega^1_Z$ is trivial. Choose a
trivialization of $\det\Omega^1_Z$, i.e., a differential form
$\omega\in H^0(Z,\det\Omega^1_Z)$ with no zeros. By
\ref{ind-smooth2}, the scheme $X:=(\cL Z)_{\red}$ is locally
nice. By \ref{dimtorsproj} and the corollary of Lemma
\ref{easy}, our trivialization of $\det\Omega^1_Z$ induces
a  trivialization $\eta$ of the dimension torsor $\Dim X$.
We put
\begin{equation}  \label{motint}
\int\limits_{Z}|\omega|:=
R\Gamma_c^{\eta}(X,\mathbb{Z}_l)\in D^b_c(\Spec k,\BZ_l),
\end{equation}
where $R\Gamma_c^{\eta}(X,\mathbb{Z}_l)$ is defined by
\eqref{RGammaceta}. Clearly $\int\limits_{Z}|\omega|$ does
not depend on the choice of $X$.

\subsection{Comparison with usual motivic integration}

In the situation of \ref{refined} (i.e., integrating a
holomorphic form with no zeros over an affinoid domain) the
usual motivic integral \cite{Lo} belongs to 
$M_k:=M'_k[L^{-1}]$, where $M'_k$ is the Grothendieck ring
of $k$-varieties\footnote{$M'_k$ is generated by elements
$[X]$ corresponding to isomorphism classes of $k$-schemes of
finite type, and the defining relations are
$[X]=[Y]+[X\setminus Y]$ for any $k$-scheme $X$ of finite
type and any closed subscheme $Y\subset X$. In particular,
these relations imply that $[X]$ depends only on the reduced
subscheme corresponding to $X$.}
and $L\in M'_k$ is the class of the affine line. Its
definition can be reformulated as follows.

Given a connected nice $k$-scheme $X$ and a trivialization
$\eta$ of its dimension torsor one chooses $\pi :X\to Y$ as in
\ref{renorm}, defines $m\in H^0(X,\BZ )=\BZ$ by
\eqref{m} and puts $[X]^{\eta}:=[Y]L^{m}\in M_k$.
If $X$ is any quasicompact quasiseparated locally nice
$k$-scheme choose closed subschemes 
$X=F_0\supset F_1\supset\ldots\supset F_n=\emptyset$
so that each $F_i$ is defined by finitely many equations
and $F_i\setminus F_{i+1}$ is nice and connected; then put
$[X]^{\eta}:=\sum_i[F_i\setminus F_{i+1}]^{\eta}$. Finally, 
in the situation of \ref{refined} one puts 
\begin{equation}  \label{motinusual}
(\int\limits_{Z}|\omega |)_{usual}:=[X]^{\eta}\in M_k,
\end{equation}
Clearly \eqref{motinusual} is well-defined, and the images of
\eqref{motinusual} and \eqref{motint}  in $K_0(D^b_c(\Spec
k,\BZ_l))$ are equal. So
\eqref{motinusual} and \eqref{motint} can be considered as
different refinements of the same object of 
$K_0(D^b_c(\Spec k,\BZ_l))$. Unless
the map $M_k\to K_0(D^b_c(\Spec k,\BZ_l))$ is
injective (which seems unlikely), the ``refined'' motivic
integral \eqref{motint}  cannot be considered as the
refinement of the usual motivic integral \eqref{motinusual}.
This is why I am using quotation marks.

\subsection{Remark}
Our definition of ``refined'' motivic integration works only
in the case of integrating a holomorphic form with no zeros
over an affinoid domain (which is probably too special for
serious applications).

On the other hand, in an unpublished manuscript
V.~Vologodsky  defined a different kind of ``refined motivic
integration'' in the case of K3 surfaces. More precisely,
let $\omega\ne 0$ be a regular differential form on a K3
surface $X$ over $F=k((t))$, $\Char k=0$. Let $A$ denote the
Grothendieck ring of the category of Grothendieck motives over
$k$, and let $I_n$ denote the motivic integral of $\omega$ 
over $X\otimes_Fk((t^{1/n}))$ viewed as an object of
$A\otimes\BQ$. Vologodsky defined objects $M_1,M_2,M_3$ of the
category of Grothendieck motives so that $I_n$ is a certain 
linear combination of the classes of $M_1,M_2,M_3$. The
objects $M_1,M_2,M_3$ depend functorially on $(X,\omega )$.
His definition of $M_1,M_2,M_3$ is mysterious.

\subsection{Counterexamples}
\label{asbadasposs}

Here are the examples promised in Remark (ii)
of \ref{affineloops}. 

\subsubsection{Not Zariski-locally trivial dimension torsor}
\label{counterex1}

Put $Y:=(\BP^1\times\BP^1)\setminus\Gamma_f$, where $\BP^1$ is
the projective line over $F:=k((t))$ and $\Gamma_f$ is the
graph of a morphism $f:\BP^1\to\BP^1$ of degree $n>0$. Clearly
$Y$ is affine, and 
\begin{equation} \label{canonbun}
\det\Omega^1_Y=p_1^*\cO(-2)\otimes
p_2^*\cO(-2)=p_1^*\cO(2n-2),
\end{equation}
where $p_1,p_2:Y\to\BP^1$ are the projections. We claim that
if $n>1$ then the dimension torsor of $\cL Y$ is not
Zariski-locally trivial. Moreover, there exists a morphism 
$\phi:\Spec R\to \cL Y$, $R:=\{ f\in k [x] |f(0)=f(1)\}$, such
that $\phi^*\Dim_{\cL Y}$ is not Zariski-locally trivial. One
constructs $\phi$ as follows. Consider the $R((t))$-module $M$
defined by \eqref{eqone}. One can represent $M$ as a direct
summand of $R((t))^2$. Indeed, the $R((t))$-module
\begin{displaymath} 
\{u=u(x,t)\in k[x]((t))^2\,|\, u(1,t)=A(t)u(0,t)\},\quad 
A(t):=\left(\begin{smallmatrix}t & 0\\ 0&
t^{-1}\end{smallmatrix}\right)
\end{displaymath}
is isomorphic to $R((t))^2$ because there exists 
$A(x,t)\in SL(2,k[x,t,t^{-1}])$ such that $A(0,t)$ is the
unit matrix and $A(1,t)=A(t)$ (to find $A(x,t)$ represent
$A(t)$ as a product of elementary matrices). Representing $M$
as a direct summand of $R((t))^2$ one gets a morphism
\begin{equation} \label{G}
g:\Spec R((t))\to\BP^1.
\end{equation}

As $p_1:Y\to\BP^1$ is a locally trivial fibration with fiber
$\BA^1$, one can represent $g$ as $p_1\varphi$ for some
$\varphi:\Spec R((t))\to Y$. Let $\phi:\Spec R\to\cL Y$ be
the morphism corresponding to $\varphi$. By \eqref{canonbun}
and the corollary of Lemma \ref{easy}, the $\BZ$-torsor
$\phi^*\Dim_{\cL Y}$ canonically identifies with the dimension
torsor of $M^{\otimes (2n-2)}$. In particular, $\phi$ has the
desired property, i.e., $\phi^*\Dim_{\cL Y}$ is not
Zariski-locally trivial. 

The class of $\phi^*\Dim_{\cL Y}$ in 
$H^1_{\et} (\Spec R,\BZ )$ is not a generator of this group
(using \eqref{canonbun} and the morphism \eqref{f} one sees
that it equals $(2n-2)v$, where $v$ is a generator). Below we
construct a slightly different pair $(Y,\phi:\Spec R\to\cL
Y)$ so that the class of $\phi^*\Dim_{\cL Y}$ in $H^1_{\et}
(\Spec R,\BZ )$ is a generator.

\subsubsection{Modification of the above example}
\label{counterex2}
Let $Y$ be the space of triples $(v,l,l')$, where $l,l'$ are
transversal 	1-dimensional subspaces in $F^2$ and $v\in l$.
Then there exists a morphism $\phi:\Spec R\to\cL Y$,
$R:=\{ f\in k [x] |f(0)=f(1)\}$, such that the class of the
$\BZ$-torsor $\phi^*\Dim_{\cL Y}$ is a generator of 
$H^1_{\et} (\Spec R,\BZ )$.

More precisely, define $\pi :Y\to\BP^1$ by $\pi (v,l,l'):=l$,
let $\tilde g:\Spec R((t))\to Y$ be such that $\pi\tilde g$
equals \eqref{G}, and let $\phi:\Spec R\to\cL Y$ be the
morphism corresponding to $\tilde g$. Then the class of 
$\phi^*\Dim_{\cL Y}$ is a generator of 
$H^1_{\et} (\Spec R,\BZ )$.

\subsubsection{Any ``unpleasant thing" can happen}
\label{unpleasanthappen}
This is what the following theorem essentially says. E.g.,
combining statement (ii) of the theorem with Weibel's
examples mentioned in \ref{Weibexamples} one sees that
for some smooth affine scheme $Y$ over $F=k((t))$ with trivial
canonical bundle there exists a reasonable closed
subscheme of $\cL Y$ which is not Zariski-locally nice (even
though its dimension torsor is trivial).

\begin{theor}  \label{impleasantoccurs}
Let $R$ be a $k$-algebra and $u\in K_{-1}(R)$.

(i) There exists a smooth affine scheme $Y$ over $F=k((t))$
and a morphism $f:\Spec R\to\cL Y$ such that the pullback of
the cotangent sheaf of $\cL Y$ to $\Spec R$ has class $u$.

(ii) If the image of $u$ in $H^1_{\et}(\Spec R,\BZ)$ equals
$0$ then one can choose $Y$ to have trivial canonical bundle
(in this case the dimension torsor of $\cL Y$ is trivial).
\end{theor}

\noindent {\bf Sketch of the proof.} Consider schemes $Y$ of
the following type\footnote{The manifold
$Y$ from \ref{counterex2} is a particular example of
\eqref{typeofmanifolds}, in which $m=n=1$ and $\dim V=1$.}:
\begin{equation} \label{typeofmanifolds}
Y=Y_0\otimes_kF,\; Y_0=(G\times V)/H,\; G=\Aut (k^m\oplus
k^n),\; H=\Aut k^m\times\Aut k^n,\; m,n\in\BN ,
\end{equation}
where $G$ and $H$ are viewed as algebraic groups over $k$ and
$V$ is a suitable representation of $H$. To prove statement
(i) of the theorem it suffices to take 
$V=\Lie [H,H]\oplus W^*$,  where $W$ is the representation of
$H$ in $k^m$. To prove (ii) it suffices to take  
$V=\Lie [H,H]\oplus W^*\oplus \det W$.

\section{Application to finite-dimensional vector
bundles on manifolds with punctures}
\label{polydisksection}

\subsection{The top cohomology}  \label{thetop}
Let $R$ be commutative,\footnote{One can formulate and prove
Theorems \ref{9.2new}-\ref{polydiskproperties} without
the commutativity assumption. In this case there is no
$S'_n$, but one can define a vector bundle on $S'_n$ to be a
collection of finitely generated projective modules $P_i$
over $\Spec R[[t_1,\ldots,t_n]][t_i^{-1}]$ with a compatible
system of isomorphisms $P_i[t_j^{-1}]\iso P_j[t_i^{-1}]$.}
$S_n:=\Spec R[[t_1,\ldots,t_n]])$, ${\bf 0}\subset S_n$ 
the subset defined by the equations $t_1=\ldots =t_n=0$,
and $S'_n :=S_n\setminus {\bf 0}$. Let $\Vect$
denote the category of vector bundles on $S'_n$ (of finite
rank).  For $L\in\Vect$ write $H^i (L)$ instead of 
$H^i(S'_n,L)$. The cohomology functors
$H^i:\Vect\to$\{$R$-modules\} vanish for
$i\ge n$ and if $n>1$ then $H^{n-1}$ commutes with base
change $R\to\tilde R$, 
$R[[t_1,\ldots,t_n]]\to\tilde R[[t_1,\ldots,t_n]]$.

\begin{theor}   \label{9.2new}
If $n>1$ then for every $\mathscr{L}\in\Vect$ the $R$-module
$H^{n-1}(S'_n,\mathscr{L})$ is \ALMOST projective. 
\end{theor}

\subsection{Derived version} \label{derivedver}
It was A.~Beilinson who explained to me that such a version
should exist. Consider
$D^\mathrm{perf}(S'_n)=\K^b(\Vect ):=$\{homotopy category of
bounded complexes in $\Vect$\}. We will decompose
$R\Gamma:\K^b(\Vect )\to D(R)$ as
$$\K^b(\Vect )\xrightarrow{(R\Gamma)_{\TTop}}\K^b(\cT_R)
\xrightarrow{\mathrm{Forget}} D(R),$$ 
where $\cT_R$ is the category of Tate $R$-modules. First of
all, we have the derived functor
\begin{equation} \label{firstfunctor}
R\Gamma:\K^b(\Vect )\to 
D^-(R[[t_1,\ldots ,t_n]])=\K^-(\cP ),
\end{equation}
where $\cP$ is the category of  projective 
$R[[t_1,\ldots ,t_n]]$-modules. Second, a projective module
$P$ over
$R[[t_1,\ldots ,t_n]]$ carries a natural topology (the
strongest one such that all
$R[[t_1,\ldots ,t_n]]$-linear maps from finitely generated
free $R[[t_1,\ldots ,t_n]]$-modules to $P$ are continuous),
so we get a functor from $\cP$ to the additive category
$R\Ttop$ of topological $R$-modules and therefore a functor
\begin{equation} \label{secondfunctor}
\K^-(\cP )\to\K^-(R\Ttop ).
\end{equation}

\begin{theor}   \label{polydiskderived}
The composition of \eqref{firstfunctor} and
\eqref{secondfunctor} belongs to the essential image of
$\K^b(\cT_R )$ in
$\K^-(R\Ttop )$, so we get a triangulated functor 
$(R\Gamma)_{\TTop}:\K^b(\Vect )\to\K^b(\cT_R)$. If 
$\cL\in\Vect$ then 
$\RGamma (\cL )\subset\K^{[0,n-1]}(\cT)$.
\end{theor}

To formulate the basic properties of $\RGamma$ we need some
notation. Let $\cC^{\Kar}$ denote the Calkin category of $R$
(see \ref{Calkincat}). Consider the functor
$\K^b(\cT_R)\to\K^b(\cC^{\Kar})$ induced by \eqref{Fi}. The
composition
\begin{equation} \label{thecomposition}
\K^b(\Vect )\xrightarrow{(R\Gamma)_{\TTop}}\K^b(\cT_R)
\to\K^b(\cC^{\Kar})
\end{equation}
will be denoted by $R\Gamma_{\discr}$, because the image of
a Tate $R$-module $T$ in $\cC^{\Kar}$ may be viewed as the
``discrete part" of $T$. One also has the ``compact part"
functor from $\cT_R$ to the category $(\cC^{\Kar})^{\circ}$
dual to $\cC^{\Kar}$: this is the composition of the
dualization functor $\cT_R\to\cT_R^{\circ}$ and the functor
$\cT_R^{\circ}\to (\cC^{\Kar})^{\circ}$ corresponding to
\eqref{Fi}. So we get 
$R\Gamma_{\comp}:\K^b(\Vect )\to\K^b((\cC^{\Kar})^{\circ})$.

Theorem \ref{9.2new} is an easy consequence of statement (i)
of the following theorem.

\begin{theor}   \label{polydiskproperties}
(i) If $\cL\in\Vect$ then $R\Gamma_{\discr} (\cL )$ is
an object of $\cC^{\Kar}$ placed in degree $n-1$.

(ii) If $\cL\in\Vect$ then $R\Gamma_{\comp} (\cL )$ is
an object of $(\cC^{\Kar})^{\circ}$ placed in degree $0$.

(iii) Let $\omega$ denote the relative (over $R$) canonical
line bundle on $S'_n$. Then there is a canonical duality
between 
$(R\Gamma)_{\TTop}(\cL)$, $\cL\in\K^b(\Vect )$, and 
$(R\Gamma)_{\TTop}(\cL^*\otimes\omega [n-1])$.
\end{theor}

\subsection{The dimension torsor corresponding to a vector
bundle}
\label{smysl}

Let $X=\Spec R$ be an affine scheme of finite type over
$\BC$.  Let $Y$ denote $X(\BC )$ equipped with the usual
topology. Given a vector bundle $\cL$ on 
$X\times (\BA^n\setminus\{ 0\} )$, $n>1$, one has the
$R$-module $M:=H^{n-1}(X\times (\BA^n\setminus\{ 0\}),\cL )$,
which is \allmost projective by Theorem \ref{9.2new}.
So by \ref{APdimtors} one has the dimension torsor $\Dim_M$
(which can be viewed as a torsor on $Y$) and its canonical
upper semicontinuous section $d^{\can}$. Here is a geometric
description of $(\Dim_M, d^{\can})$.

(i) Notice that a complex vector bundle of any rank $m$ on 
the topological space $\BC^n\setminus\{ 0\}$ is trivial
(because $\pi_{2n-2}(GL(m,\BC ))=0$), and the homotopy
classes of its trivializations form a torsor over
$\pi_{2n-1}(GL(m,\BC ))$. One has the natural morphism
$\pi_{2n-1}(GL(m,\BC ))\to \pi_{2n-1}(GL(\infty ,\BC ))=
K_0^{\Top}(S^{2n-2})=\BZ$. So $\cL$ defines a
$\BZ$-torsor $T_{\cL}$ on $Y$. 

(ii) More generally, a finite
complex $\cL^{\bcdot}$ of topological vector bundles on
$Y\times (\BC^n\setminus\{ 0\} )$ defines a $\BZ$-torsor
$T_{\cL^{\bcdot}}:=\sum_i(-1)^iT_{\cL^i}$, and a homotopy
equivalence $f:\cL_1^{\bcdot}\to\cL_2^{\cdot}$ defines an
isomorphism $T_{\cL_1}\iso T_{\cL_2}$ (because the dimension
torsor of $\Cone (f)$ is canonically trivialized). Of course,
this isomorphism depends only on the homotopy class of $f$.
An extension of $\cL^{\cdot}$ to an object of the
homotopy category of complexes of topological vector bundles
on $Y\times \BC^n$ defines a
trivialization of $T_{\cL^{\bcdot}}$.

(iii) Let $\cL$ be an algebraic vector bundle on 
$X\times (\BA^n\setminus\{ 0\} )$, $n>1$. Let
$j:\BA^n\setminus\{ 0\}\to\BA^n$ be the embedding. For each
$x\in X$ the sheaf $j_*\cL_x$ is coherent
and has a finite locally free resolution (here $\cL_x$ is
the restriction  of $\cL$ to 
$\{ x\}\times (\BA^n\setminus\{ 0\} )$). So
by (ii) one gets a trivialization of $T_{\cL_x}$ for each
$x$, i.e., a set-theoretical section $s$ of $T_{\cL}$.

(iv) One can show that $(T_{\cL},s)$ is canonically
isomorphic to $(\Dim_M, d^{\can})$ (maybe up to a sign).

\subsection{Central extension}   \label{centralext}

\begin{theor}  \label{apcoh}
Let $X$ be a smooth scheme over $S:=\Spec R$ of pure relative
dimension $n>1$. Let $F\subset X$ be a closed subscheme which
is finite  over $S$ and a locally complete intersection over
$S$. Let $j:X\setminus F\mono X$ denote the open embedding.
Then for any vector bundle $\cL$ on $X\setminus F$ the
$R$-module $H^0(X, \RR^{n-1}j_*\cL)$ is \ALMOST projective.
\end{theor}

This easily follows from Theorems \ref{9.2new} and \ref{3}(i).

Now let $O_F$ be the ring of regular functions on the
formal completion of $X$ along $F$. In the situation of
Theorem \ref{apcoh} 
$H^0(X, \RR^{n-1}j_*\cL)$ is an $O_F$-module, so it is
equipped with an action of the group scheme $G:=O_F^{\times}$.
Therefore applying
\eqref{central2} one gets a central extension
\begin{equation}    \label{central3}
0\to\BG_m\to\widehat{G}_{\cL}\to O_F^{\times}\to 0.
\end{equation}

\noindent {\bf Remarks.} (i) Suppose that in the situation of
Theorem \ref{apcoh} $\cL$ extends to a vector bundle on $X$.
Then the $R$-module $H^0(X, \RR^{n-1}j_*\cL)$ is projective,
and therefore the extension \eqref{central3} canonically
splits.

(ii) Suppose that in the situation
of Theorem \ref{apcoh} $F\subset\tilde F\subset X$ and
$\tilde F$ satisfies the same conditions as $F$. Put
$\tilde\cL:=\cL|_{X\setminus\tilde F}$. Then we have the
central extension \eqref{central3} and a similar central
extension
\begin{equation}    \label{central4}
0\to\BG_m\to\widehat{G}_{\tilde\cL}\to 
O_{\tilde F}^{\times}\to 0.
\end{equation}
Using the functor $(R\Gamma)_{\TTop}$ from \ref{derivedver}
one can construct a canonical morphism from \eqref{central4}
to \eqref{central3} which induces the restriction map
$O_{\tilde F}^{\times}\to O_F^{\times}$ and the identity map
$\BG_m\to\BG_m$. If $F=\emptyset$ this amounts to Remark (i)
above.
 
\subsection{Commutativity theorem}
Let $c_{\cL}:O_F^{\times}\times O_F^{\times}\to\BG_m$ be the
commutator map of the central extension \eqref{central3}.

\begin{theor}        \label{commutativity}
Suppose that in the situation of Theorem \ref{apcoh} $n=2$.
Then $c_{\cL}=1$ if and only if $\det\cL$ extends to an
invertible sheaf on $X$.
\end{theor}

\noindent {\bf Remarks.} (i) If an extension of $\det\cL$ to
an invertible sheaf on $X$ exists it equals
$j_*\det\cL$. In particular, the extension is unique.

(ii) Theorems \ref{apcoh} and \ref{commutativity} are
still true for vector bundles on $(\Spec O_F)\setminus F$
instead of $X\setminus F$. 

\medskip

\noindent {\bf Question.} What is the geometric meaning of
$c_{\cL}$ and the condition $c_{\cL}=1$ if $n>2$?

\subsection{Generalizing the notion of vector bundle on a
surface}  \label{generalizing}
Let $G$ be a reductive group over $\BQ$. The moduli scheme of
$G$-bundles on $\BP^2_{\BQ}$ trivialized over a fixed
projective line $\BP^1_{\BQ}\subset \BP^2_{\BQ}$ has a
remarkable ``Uhlenbeck compactification'' $\fU_G$ constructed
in \cite{FGK,BFG}, which goes back to the physical picture of
``instanton gas''. It would be very important to interpret
$\fU_G$ as a moduli scheme of some kind of  geometric
objects.\footnote{E.g., such an interpretation would
hopefully allow to define an analog of $\fU_G$ for {\it
any\,} proper smooth surface.} These conjectural new objects
are, so to say, ``$G$-bundles with singularities". I suggest
to call them {\it
$G$-gundles.} The new word ``gundle'' can be considered as an
abbreviation for ``generalized $G$-bundle". On the other
hand, its first 3 letters  are also the first letters of the
names of D.~Gaitsgory, V.~Ginzburg, K.~Uhlenbeck, and
H.~Nakajima.\footnote{Gaitsgory is an author of
\cite{FGK,BFG}, and the relation of the other 3
mathematicians to these articles is explained in the
introductions to them.}

It turns out that the central extension \eqref{central3}
allows to give a definition of $GL(n)$-gundle on any
smooth family of surfaces over any base so that $\fU_{GL(n)}$
identifies with the moduli scheme of $GL(n)$-gundles on
$\BP^2_{\BQ}$ trivialized over $\BP^1_{\BQ}$.

Let $X$ be a scheme smooth over $S$ of pure relative
dimension 2. The definition of $GL(n)$-gundle on $X$ consists
of several steps. I will only explain the first one and list
the other steps.

\subsubsection{Pre-gundles 1} 
Let $F$ be as in Theorem \ref{apcoh}. 

\medskip

\noindent {\bf Definition.} A {\it  $GL(n)$-pre-gundle on
$X$ nonsingular outside $F$\,} is a pair that consists of
a rank $n$ vector bundle $\cL$ on $X\setminus F$ and a
splitting of $\eqref{central3}$. The groupoid of such pairs
will be denoted by $\Pregun_F(X)$.

\medskip

\noindent {\bf Remark.} If $\eqref{central3}$ admits a
splitting then by Theorem \ref{commutativity} $\det\cL$
extends to a line bundle on $X$.

\subsubsection{Pre-gundles 2}
If $F$, $\tilde F$ are as in Theorem \ref{apcoh} and 
$\tilde F\supset F$ then one defines a fully faithful functor
$\Pregun_F(X)\to\Pregun_{\tilde F}(X)$ using Remark (ii) at
the end of \ref{centralext}.

\subsubsection{Pre-gundles 3} If $X$ is projective over $S$
one defines the groupoid of pre-gundles on $X$ to be the
inductive 2-limit of $\Pregun_F(X)$ over all closed subschemes
$F\subset X$ which are finite over $S$ and locally complete
intersections over $S$. This groupoid is denoted by
$\Pregun (X)$, and its objects are called 
{\it $GL(n)$-pre-gundles on $X$.}

If $X$ is arbitrary one first defines $\Pregun_F(X)$ for
{\it any\,} subscheme $F\subset X$ quasi-finite over $S$ (a
standard etale or Nisnevich localization technique allows to
reduce this to the case of finite locally complete
intersection considered above). Then one defines 
$\Pregun (X)$ to be the inductive 2-limit of $\Pregun_F(X)$
over {\it all\,} closed subschemes
$F\subset X$ quasi-finite over $S$.

\subsubsection{Remark}    \label{fieldcase}
If $S$ is the spectrum of a field then $\Pregun (X)$
identifies with the groupoid of pairs $(\cL ,Z)$ with $\cL$
being a $GL(n)$-bundle on $X$ and $Z$ a $0$-cycle on $X$,
it being understood that an isomorphism 
$(\cL_1 ,Z)\iso (\cL_2 ,Z)$ is the same as an isomorphism
$\cL_1\iso\cL_2$ and that there are no isomorphisms 
$(\cL_1 ,Z_1)\iso (\cL_2 ,Z_2)$ if $Z_1\ne Z_2$.

To see this, first notice that for any finite $F\subset X$ a
vector bundle on $X\setminus F$ uniquely extends to $X$.
Second, by Remark (i) from \ref{centralext}, the central
extension \eqref{central3} has a canonical splitting, so
{\it all\,} splittings of \eqref{central3} are parametrized
by $\Hom (O_F^{\times},\BG_m )$, i.e., by the group of
$0$-cycles on $X$ supported on $F$.

\subsubsection{Pre-gundles 4} Let $S$ be again arbitrary.
Associating to an $S$-scheme $S'$ the groupoid of pre-gundles
on $X\times_SS'$ one gets a (non-algebraic) $S$-stack
$\Pregun_X$. This is {\it the stack of $GL(n)$-pre-gundles on
$X$.}

\subsubsection{Gundles} One defines a closed substack
$\Gun_X\subset\Pregun_X$, whose formation commutes with base
change $S'\to S$. Its $S$-points are called
{\it $GL(n)$-gundles on $X$.}

By \ref{fieldcase}, if $S=\Spec k$ with $k$ being a field then
$GL(n)$-pre-gundles on $X$ identify with pairs $(\cL ,Z)$ with
$\cL$ being a $GL(n)$-bundle on $X$ and $Z$ being a $0$-cycle
on $X$. It turns out that {\it such a pair $(\cL ,Z)$ is a
$GL(n)$-gundle if and only if $Z\ge 0$.}

\medskip

\noindent {\bf Remark.} I can define the closed substack
$\Gun_X\subset\Pregun_X$ using the method of \cite{FGK,BFG},
i.e., by working with various curves on $X$. Unfortunately, I
do not know a ``purely 2-dimensional" way to do it.

\subsubsection{Hope} If $X$ is proper over $S$ then
$\Gun_X$ is an algebraic stack.

\subsubsection{Fact} Now let $S=\Spec\BQ$ and
$X=\BP^2_{\BQ}$. Fix a projective line
$\BP^1_{\BQ}\subset\BP^2_{\BQ}$, and consider the open
substack $U\subset\Gun_X$ parametrizing those gundles which
are nonsingular on a neighborhood of $\BP^1$ and whose
restriction to $\BP^1$ is trivial. Then
$U$ identifies with the quotient of the ``Uhlenbeck
compactification'' $\fU_{GL(n)}$ from \cite{BFG} by the
action\footnote{$GL(n)$ acts on $\fU_{GL(n)}$ by changing the
trivialization over $\BP^1$.} of $GL(n)$. In particular, the
stack $U$ is algebraic.


\bibliographystyle{alpha}

\end{document}